\documentclass[11pt,a4paper]{article}
\usepackage[margin=28mm]{geometry}
\usepackage{amsmath,amssymb,amsthm,mathtools,mathrsfs}
\usepackage[expansion=false]{microtype}
\usepackage[hidelinks]{hyperref}
\allowdisplaybreaks
\numberwithin{equation}{section}
\newtheorem{theorem}{Theorem}[section]
\newtheorem{lemma}{Lemma}[section]
\newtheorem{proposition}[lemma]{Proposition}
\newtheorem{corollary}[lemma]{Corollary}
\theoremstyle{remark}
\newtheorem{remark}[lemma]{Remark}
\newcommand{\R}{\mathbb R}
\newcommand{\TT}{\mathbb T^2}
\newcommand{\dd}{\,\mathrm dx}
\newcommand{\diver}{\operatorname{div}}
\newcommand{\curl}{\operatorname{curl}}
\newcommand{\Dt}{D_t}
\newcommand{\eps}{\varepsilon}
\title{Global strong solutions and large-time behavior
for the two-dimensional Navier--Stokes--Poisson system
with general bulk viscosity}
\author{Jing Han$^{a,b}$, Yanfang Peng$^a$ \thanks{Email: jinghanfc@163.com(J. Han), pyfang2005@sina.com(Y. Peng)} \\[3mm] \small \emph{ a. School of Mathematical Sciences, Guizhou Normal University, Guiyang 550025, P. R. China} \\ \small \emph {b.  School of Mathematics and Physics, Anshun University, Anshun 561000, P. R. China}}
\date{}
\begin{document}
\maketitle
\begin{abstract}
The global existence and uniqueness of strong solutions to the
two-dimensional repulsive Navier--Stokes--Poisson system with general
density-dependent bulk viscosity is established for arbitrarily large
initial data allowing vacuum. The periodic problem and the Cauchy
problem with a positive constant background are treated under a common
growth condition on the viscosity primitive. The coefficient class
includes every superlinear power law and near-linear laws with
logarithmic gains. A density upper bound and low-order estimates,
uniform in time and independent of a positive density lower bound,
are obtained without an initial acceleration compatibility condition.
Exponential relaxation is proved on the torus, and qualitative
convergence to equilibrium is established in the whole space under a
finite electric energy condition and without spatial moment assumptions.
If the initial density vanishes at a point, the vacuum persists and
density gradients become unbounded as time tends to infinity.
\end{abstract}



\section{Introduction and main results}\label{pc:introduction}

The purpose of this paper is to establish global strong solvability
and large-time behavior for the two-dimensional repulsive
Navier--Stokes--Poisson (NSP) system with general bulk viscosity.
The initial data may have arbitrary amplitude and may contain vacuum.
The main issue is to control the density uniformly in time while
retaining the electric energy and the low-frequency part of the flow.
Let $\Omega=\TT:=\R^2/\mathbb Z^2$ or $\Omega=\R^2$. We consider
\begin{equation}\label{pc:model}
\begin{cases}
 \rho_t+\diver(\rho u)=0,\\
 \rho\dot u+\nabla\rho^\gamma
 =\mu\Delta u+\nabla((\mu+\lambda(\rho))\diver u)+\rho\nabla\Phi,\\
 \Delta\Phi=\rho-1.
\end{cases}
\end{equation}
Here $\rho$, $u$, and $\Phi$ denote the density, velocity, and
electrostatic potential. The shear viscosity $\mu>0$ and pressure
exponent $\gamma>1$ are fixed. We write
\[
 P=\rho^\gamma,\qquad E=\nabla\Phi,\qquad
 \Dt=\partial_t+u\cdot\nabla,\qquad \dot u=\Dt u.
\]
The initial density and momentum are prescribed by
\begin{equation}\label{pc:initial-traces}
 \rho|_{t=0}=\rho_0,\qquad (\rho u)|_{t=0}=\rho_0u_0.
\end{equation}
On the unit torus, the total mass is one and $\bar\Phi=0$.
In the whole space, the background and far-field densities are both
one, and the electric field has finite energy. The signs in
\eqref{pc:model} give the repulsive coupling and a positive electric
contribution to the energy. Spatial integrals and
norms are over the domain under consideration; on $\TT$ we set
$\bar f=\int_{\TT}f\dd$. Our curl convention is
\[
 \omega=\curl u=\partial_2u_1-\partial_1u_2,\qquad
 \nabla^\perp=(\partial_2,-\partial_1).
\]

For the two-dimensional barotropic Navier--Stokes equations without
an electric force, density-dependent bulk viscosity provides a mechanism
for controlling large solutions. For $\lambda(\rho)=\rho^\beta$,
Vaigant and Kazhikhov~\cite{VaigantKazhikhov1995} proved global
solvability for large positive periodic data when $\beta>3$.
Huang and Li~\cite{HuangLi2016} obtained global strong solutions
allowing vacuum for $\beta>4/3$.
The general-viscosity theory of Fan, Jiang, Li, and Lv~\cite{SourceNS}
replaces the power-law restriction by a growth condition on the viscosity
primitive. These results concern the uncoupled Navier--Stokes equations.
The Poisson force in \eqref{pc:model} changes both the energy estimates
and the large-time analysis.

For the repulsive Navier--Stokes--Poisson system, the electric field
also affects the low-frequency dynamics. Li, Matsumura, and
Zhang~\cite{LiMatsumuraZhang2010} studied this effect and optimal decay
for small perturbations in three dimensions. Chikami and
Danchin~\cite{ChikamiDanchin} and Shou and Song~\cite{ShouSong}
established global critical-space theories near a positive constant
equilibrium. In two-dimensional bounded domains, Du and
Zhao~\cite{DuZhao2026} proved large-data strong solvability allowing
vacuum for power-law bulk viscosity with $\beta>4/3$ and a positive,
possibly nonconstant background. The present paper treats the periodic
and whole-space problems with general bulk viscosity. In particular,
the whole-space analysis uses finite relative energy without a spatial
moment condition. Precise comparisons of the models and data classes
are given after the main results.

\subsection{Bulk viscosity}

Put $h(r)=2\mu+\lambda(r)$. The bulk viscosity satisfies
\begin{equation}\label{pc:lambda}
 \lambda\in C^1([0,\infty))\cap C^2((0,\infty)),\quad
 \lambda\ge0,\quad r|\lambda'(r)|\le C_\lambda h(r)
 \quad(r>0).
\end{equation}
Define its primitive by
\[
 \theta(r)=\int_1^r\frac{h(z)}z\,\mathrm dz,\qquad r>0.
\]
For some $\kappa\in(1/2,1]$, $c_\lambda>0$, and $r_\lambda\ge2$,
we assume
\begin{equation}\label{pc:growth}
 \theta(r)\ge c_\lambda r[\log(e+r)]^\kappa,
 \qquad r\ge r_\lambda.
\end{equation}
By~\cite[Lemma~2.2]{SourceNS}, with $\ell=\max\{1,C_\lambda\}$,
\begin{equation}\label{pc:coefficient}
 h(r)\le C(1+r^\ell),\qquad
 \lambda(r)\le C(1+\theta(r)_+).
\end{equation}
The relative internal energy is
\begin{equation}\label{pc:internal}
 H(r\mid1)=\frac{r^\gamma-1-\gamma(r-1)}{\gamma-1}.
\end{equation}

The hypotheses do not require monotonicity of $\lambda$ or
$\lambda(0)=0$. In particular, they include
$\lambda(r)=b r^\beta$ for every $b>0$ and $\beta>1$, and
$\lambda(r)=b r[\log(e+r)]^\eta$ for every $b>0$ and $\eta>1/2$.
Their common role is the lower growth bound for $\theta$, rather
than a fixed power-law formula for $\lambda$.

\subsection{The periodic problem}

Our first result gives global solvability and exponential relaxation
toward the equilibrium determined by the total mass and momentum.

\begin{theorem}[Periodic strong solution and relaxation]\label{pcstr:main}
Let $q>2$. Assume \eqref{pc:lambda}--\eqref{pc:growth}, and let
\[
 0\le\rho_0\in W^{1,q}(\TT),\qquad
 \int_{\TT}\rho_0\dd=1,\qquad u_0\in H^1(\TT).
\]
Then the periodic system \eqref{pc:model}, with initial density
$\rho_0$ and initial momentum $\rho_0u_0$, has a unique global strong
solution $(\rho,u,\Phi)$. The potential is normalized by $\bar\Phi=0$ and is
determined by the density. For every finite $T>0$ and
$1<s<q/(q-1)$,
\begin{equation}\label{pcstr:strong-class}
\begin{aligned}
 &\rho\in C([0,T];W^{1,q}),\qquad
            \rho_t\in L^\infty(0,T;L^2),\\
 &u\in L^\infty(0,T;H^1)\cap L^2(0,T;H^2)
                         \cap L^s(0,T;W^{2,q}),\\
 &\sqrt t\,u\in L^2(0,T;W^{2,q}),\qquad
   \sqrt t\,u_t\in L^2(0,T;H^1),\qquad
   \sqrt\rho\,u_t\in L^2(0,T;L^2),\\
 &\rho u,\sqrt\rho\,u\in C([0,T];L^2),\\
 &\nabla\Phi\in C([0,T];W^{2,q}),\qquad
   \nabla\Phi_t\in L^\infty(0,T;H^1).
\end{aligned}
\end{equation}
No initial acceleration compatibility condition is imposed. Moreover,
\begin{equation}\label{pcstr:uniform-bound}
 \sup_{t\ge0}\left(\|\rho(t)\|_{L^\infty}+\|u(t)\|_{H^1}\right)\le C.
\end{equation}
The conserved velocity is
\[
 u_\infty=\int_{\TT}\rho_0u_0\dd.
\]
For every $1\le p<\infty$,
\begin{equation}\label{pcstr:decay}
 \|\rho(t)-1\|_{L^p}+\|u(t)-u_\infty\|_{W^{1,p}}
                 +\|\nabla\Phi(t)\|_{W^{1,p}}
       \le C_pe^{-c_pt},\qquad t\ge1.
\end{equation}
The constant in \eqref{pcstr:uniform-bound} depends only on
$\|\rho_0\|_{L^\infty}$, $\|u_0\|_{H^1}$, $\mu$, $\gamma$, the
fixed torus, the constants in \eqref{pc:lambda}--\eqref{pc:growth},
and a bound for $\lambda$ on $[0,2]$. The decay constants also
depend on $p$. They are independent of a positive lower density
bound and of $\|\rho_0\|_{W^{1,q}}$.
\end{theorem}

\begin{remark}[Comparison with the Navier--Stokes equations]
\label{pc:remark-viscosity}
The following three works concern the barotropic Navier--Stokes
equations without a Poisson force.
Fan, Jiang, Li, and Lv~\cite{SourceNS} treat general bulk viscosity
on $\TT$ and on $\R^2$, with far-field vacuum in the Cauchy problem.
Their bounded-domain paper~\cite{SourceSlipNS} treats
$\lambda(\rho)=\rho^\beta$, $\beta>1$, under
$u\cdot n=0$ and $\curl u=0$, and proves uniform bounds,
exponential relaxation, and vacuum-induced gradient growth.
Jiu, Wang, and Xin~\cite{JiuWangXinPositive} consider
positive far-field density on $\R^2$, with power-law bulk viscosity
and weighted energy data.
Here the equation is instead the coupled NSP system
\eqref{pc:model}; its electric energy and force require additional
estimates.
\end{remark}

\subsection{The Cauchy problem}

Let $q>2$. We assume
\begin{equation}\label{pcws:data}
 \rho_0\ge0,\qquad
 \rho_0-1\in H^1(\R^2)\cap W^{1,q}(\R^2)\cap\dot H^{-1}(\R^2),
 \qquad u_0\in H^1(\R^2).
\end{equation}
Here
$\|f\|_{\dot H^{-1}}=\||\xi|^{-1}\widehat f\|_{L^2}$
for the unitary Fourier transform. Thus the last density condition
is precisely the finite-energy condition for
$E_0=\nabla\Delta^{-1}(\rho_0-1)$. No smallness is imposed.

\begin{theorem}[Global strong solutions on $\R^2$]\label{pcwstr:main}
Assume \eqref{pc:lambda}--\eqref{pc:growth} and let the initial data
satisfy \eqref{pcws:data}. Then \eqref{pc:model} admits a unique
global strong solution $(\rho,u,E)$ with initial density $\rho_0$, initial momentum
$\rho_0u_0$, and initial field
$E_0=\nabla\Delta^{-1}(\rho_0-1)$. For every finite $T>0$ and
$1<s<q/(q-1)$,
\begin{equation}\label{pcwstr:class}
\begin{aligned}
 &\rho-1\in C([0,T];H^1\cap W^{1,q}),\qquad
              \rho_t\in L^\infty(0,T;L^2),\\
 &u\in C((0,T];H^1)\cap L^\infty(0,T;H^1),\\
 &u\in L^2(0,T;H^2)\cap L^s(0,T;W^{2,q}),\\
 &\sqrt t\,u\in L^2(0,T;W^{2,q}),\qquad
   \sqrt t\,u_t\in L^2(0,T;H^1),\qquad
   \sqrt\rho\,u_t\in L^2(0,T;L^2),\\
 &\rho u,\sqrt\rho\,u\in C([0,T];L^2),\\
 &E\in C([0,T];H^2\cap W^{2,q}),\qquad
   E_t\in L^\infty(0,T;H^1).
\end{aligned}
\end{equation}
The field is the unique curl-free $L^2$ field with
$\diver E=\rho-1$. The solution satisfies the relative energy
identity \eqref{pcws:energy} and the time-independent bound
\begin{equation}\label{pcwstr:uniform-bound}
 \sup_{t\ge0}\left(\|\rho(t)\|_{L^\infty}
                         +\|u(t)\|_{H^1}\right)\le C,
\end{equation}
and the modified acceleration satisfies
\begin{equation}\label{pcwstr:modified-budget}
 \int_0^\infty\|\sqrt\rho(\dot u-E)\|_{L^2}^2\,\mathrm dt\le C.
\end{equation}
Moreover, as $t\to\infty$,
\begin{equation}\label{pcwstr:relaxation}
\begin{aligned}
 &\|\rho(t)-1\|_{L^p}+\|\nabla u(t)\|_{L^p}\longrightarrow0,
                           &&2\le p<\infty,\\
 &\|u(t)\|_{W^{1,p}}+\|E(t)\|_{W^{1,p}}\longrightarrow0,
                           &&2<p<\infty.
\end{aligned}
\end{equation}
Arbitrary initial energy and vacuum are allowed, and no initial
acceleration compatibility condition is required. The constants in
\eqref{pcwstr:uniform-bound}--\eqref{pcwstr:modified-budget} depend on
the initial relative energy, $\|\rho_0\|_{L^\infty}$,
$\|u_0\|_{H^1}$, and the fixed viscosity and pressure data;
they are independent of a positive density lower bound.
\end{theorem}

\begin{remark}[Perturbative NSP theory]\label{pc:remark-large-data}
Chikami and Danchin~\cite{ChikamiDanchin} consider the repulsive
NSP Cauchy problem on $\R^d$, $d\ge2$, with a positive
constant background. Shou and Song~\cite[Theorem~1.1]{ShouSong}
also consider this repulsive coupling, allowing density-dependent
shear and bulk viscosities. Their viscous operator is
\[
 \diver\bigl(\mu_1(\rho)(\nabla u+\nabla u^{\mathsf T})\bigr)
       +\nabla\bigl(\mu_2(\rho)\diver u\bigr).
\]
Taking $\mu_1\equiv\mu$, $\mu_2=\lambda$, and $P=\rho^\gamma$,
and normalizing the background and coupling strength to one,
gives the equation form in \eqref{pc:model}.
Their global critical-space results require small perturbations
of a positive equilibrium. Theorem~\ref{pcwstr:main} concerns a
different regime: two dimensions, constant shear viscosity,
bulk viscosity satisfying \eqref{pc:lambda}--\eqref{pc:growth},
and large data allowing vacuum. The data and coefficient classes
therefore differ.
\end{remark}

\begin{remark}[The bounded-domain NSP problem]\label{pc:remark-bounded-nsp}
Du and Zhao~\cite[Theorem~1.2]{DuZhao2026} study the repulsive
NSP system in a two-dimensional bounded simply connected domain,
with Navier-slip velocity conditions, a Neumann potential, and a
positive possibly nonconstant background. They prove global strong
solvability with large data and vacuum for
$\lambda(\rho)=\rho^\beta$, $\beta>4/3$.
Our theorems concern $\TT$ and $\R^2$ with constant background one.
In these settings the viscosity assumptions include every
$\beta>1$, near-linear laws
$\lambda(r)=r[\log(e+r)]^\eta$, $\eta>1/2$, and nonmonotone
coefficients satisfying \eqref{pc:lambda}--\eqref{pc:growth}.
\end{remark}

\begin{remark}[The whole-space NSP estimates]\label{pc:remark-method}
For the linearization of \eqref{pc:model} on $\R^2$,
Proposition~\ref{pcws:linear-acceleration} shows that finite electric
energy does not imply $u_t\in L^2(0,\infty;L^2)$.
The modified acceleration $\dot u-E$, also used for the bounded-domain
NSP problem in~\cite{DuZhao2026}, yields
\eqref{pcwstr:modified-budget}.
The additional whole-space estimates are uniform density moments
and relative-pressure dissipation obtained before the density
maximum is bounded, with no spatial moment condition in
\eqref{pcws:data}. Bessel operators enter only these auxiliary
estimates; the Poisson equation remains unchanged.
\end{remark}

\begin{remark}[Relaxation and vacuum for the present NSP system]
\label{pc:remark-relaxation}
For \eqref{pc:model}, Theorem~\ref{pcstr:main} gives exponential
relaxation on $\TT$, while Theorem~\ref{pcwstr:main} gives
qualitative convergence on $\R^2$ with far-field density one.
No whole-space decay rate or $L^2$ convergence of $u$ and $E$ is
asserted. Corollaries~\ref{pcstr:vacuum} and~\ref{pcwr:vacuum}
show that an initial vacuum point persists. Thus density
relaxation in the stated finite-$L^p$ ranges coexists with
$\|\rho(t)-1\|_{L^\infty}\ge1$ and unbounded density gradients.
\end{remark}

\begin{remark}[Further questions for the NSP system]\label{pc:remark-future}
For \eqref{pc:model} on $\R^2$, natural further questions are
quantitative decay under additional low-frequency assumptions and
$L^2$ relaxation of $u$ and $E$ in the present data class.
The endpoint $\kappa=1/2$ is outside the current moment argument,
which requires $(2\kappa-1)k\ge1$; this does not establish a sharp
solvability threshold. Extending the present periodic and whole-space
estimates to nonconstant positive backgrounds would require
relative energy and flux estimates around a nonconstant stationary
state.
\end{remark}

\subsection{Outline of the proof}

The common step is a time-independent upper bound for the density.
We use the viscosity primitive and the abstract level argument
of~\cite{SourceNS}, after establishing the required moment, pressure,
and effective-flux estimates for the coupled system. The electric
energy is included throughout. On the torus, conservation of momentum
and the Poincar\'e inequality provide velocity control, and a relative
energy argument yields exponential relaxation.

On $\R^2$, the periodic velocity-potential bound is unavailable,
and even the linearized system need not have a finite global time
integral of the ordinary acceleration under the electric energy
assumption. Appendix~\ref{pc:low-frequency} records these obstructions.
A screened auxiliary potential yields density moments, and a separate
argument gives relative-pressure dissipation, both before the density
maximum is controlled. Estimates for the modified acceleration
$\dot u-E$ then provide the effective-flux cost needed for the density
bound. Time-weighted derivative estimates then permit
positive approximation and passage to vacuum without an initial
acceleration compatibility condition. The resulting dissipation gives
qualitative whole-space relaxation.

Sections~\ref{pc:periodic-estimates} and~\ref{pc:periodic-strong}
prove the periodic result. Sections~\ref{pc:whole-estimates}
and~\ref{pc:whole-strong} treat the whole-space problem.

\section{Periodic a priori estimates}\label{pc:periodic-estimates}
\subsection{Energy and elliptic estimates}

We first prove the uniform estimates needed for
Theorem~\ref{pcstr:main}. All integrals in this section are over the
unit torus. Let a smooth positive solution exist on $[0,S]$, with
\[
 \int\rho_0\dd=1,\qquad
 \|\rho_0\|_{L^\infty}+\|u_0\|_{H^1}\le C_0.
\]
The constants below may depend on $C_0,\mu,\gamma$, the fixed torus,
the constants in \eqref{pc:lambda}--\eqref{pc:growth}, and a bound for
$\lambda$ on $[0,2]$. They are independent of $S$ and of the positive
lower density bound. Constants with a fixed exponent as subscript may
also depend on that exponent.

The electric force has zero integral. In fact, periodic integration by
parts gives, componentwise,
\begin{equation}\label{pcper:electric-mean}
 \int\rho\nabla\Phi\dd
 =\int(1+\Delta\Phi)\nabla\Phi\dd=0,
 \qquad \frac{\mathrm d}{\mathrm dt}\int\rho u\dd=0.
\end{equation}
Put $u_\infty=\int\rho_0u_0\dd$. The change of variables
\[
 (\rho,u,\Phi)(t,x)\longmapsto
 (\rho(t,x+tu_\infty),u(t,x+tu_\infty)-u_\infty,
                         \Phi(t,x+tu_\infty))
\]
preserves \eqref{pc:model} and $\bar\Phi=0$. We use this frame throughout
the estimates, so that $\int\rho u\dd=0$ and $|u_\infty|\le C$.
The material acceleration and all spatial derivatives are unchanged
by this transformation.

Set
\begin{equation}\label{pcper:definitions}
\begin{gathered}
 G=h(\rho)\diver u-(P-1),\qquad
 R_t=\max\{2,\sup_{0\le s\le t}\|\rho(s)\|_{L^\infty}\},\\
 A_1^2=\int\left(\mu\omega^2+\frac{G^2}{h(\rho)}\right)\dd,
 \qquad A_2^2=\int\rho|\dot u|^2\dd,\\
 A_3^2=\int\left(h(\rho)(\diver u)^2+\mu\omega^2\right)\dd,
 \qquad \mathscr E=e+A_1^2.
\end{gathered}
\end{equation}
The flux $G$ is centered at the equilibrium pressure, so its energy
identity contains no mean-pressure time derivative.

\begin{lemma}[Basic energy and weighted means]\label{pcper:basic-lemma}
One has
\begin{equation}\label{pcper:energy}
 \mathcal H=\int\left(\frac12\rho|u|^2+H(\rho\mid1)
                                  +\frac12|\nabla\Phi|^2\right)\dd,
 \qquad \mathcal H'+A_3^2=0,
\end{equation}
and
\begin{equation}\label{pcper:coercivity}
 \sup_{t\le S}\int(\rho|u|^2+\rho^\gamma+|\nabla\Phi|^2)\dd
       +\int_0^S A_3^2\,\mathrm dt\le C,
 \qquad \|u\|_{H^1}\le CA_3.
\end{equation}
Moreover,
\begin{equation}\label{pcper:acceleration-mean}
 \int\rho\dot u\dd=0,\qquad
 \|\dot u\|_{L^m}\le C\sqrt m\,\|\nabla\dot u\|_{L^2},
 \quad m\ge2.
\end{equation}
\end{lemma}
\begin{proof}
The Poisson and continuity equations imply
\[
 \Delta\Phi_t=-\diver(\rho u),\qquad
 \int\rho u\cdot\nabla\Phi\dd
       =-\frac12\frac{\mathrm d}{\mathrm dt}\|\nabla\Phi\|_{L^2}^2.
\]
Testing momentum by $u$ therefore gives \eqref{pcper:energy}.
The mass condition makes the relative and absolute internal energies
differ by a fixed constant. The initial electric energy is bounded
by $C\|\rho_0-1\|_{L^2}^2$. This proves the first estimate in
\eqref{pcper:coercivity}. Also
\[
 |\bar u|=\left|\int\rho(u-\bar u)\dd\right|
 \le\|\rho\|_{L^\gamma}
       \|u-\bar u\|_{L^{\gamma/(\gamma-1)}}\le C\|\nabla u\|_{L^2}.
\]
Use $\int|\nabla u|^2=\int((\diver u)^2+\omega^2)$ to obtain
the remaining estimate. Integrating momentum and using
\eqref{pcper:electric-mean} gives $\int\rho\dot u=0$.
The same mean argument and the two-dimensional Sobolev estimate with
constant $C\sqrt m$ apply to $\dot u$.
\end{proof}

For a vector field $f$, let $Qf=(-\Delta)^{-1}\diver f$ have zero
mean. Then
\begin{equation}\label{pcper:potential-estimates}
\begin{gathered}
 Q\nabla g=-g+\bar g,\qquad
 \|\nabla Qf\|_{L^p}\le Cp\|f\|_{L^p}\quad(p\ge2),\\
 \|Qf\|_{L^{2r/(2-r)}}\le C_r\|f\|_{L^r}\quad(1<r<2).
\end{gathered}
\end{equation}
These follow from periodic Calder\'on--Zygmund and Sobolev estimates;
weak $(1,1)$, interpolation, and duality give the displayed $Cp$ bound.
For zero-mean scalar functions we also use
\begin{equation}\label{pcper:GN}
\begin{gathered}
 \|f\|_{L^p}\le C\sqrt p\,\|f\|_{L^2}^{2/p}
                       \|\nabla f\|_{L^2}^{1-2/p},\quad p\ge2,
 \\
 \int\exp\left(\frac{c|v|^2}{\|v\|_{H^1}^2}\right)\dd\le C.
\end{gathered}
\end{equation}
The second bound, for $v\ne0$, follows by expanding the exponential
and using $\|v\|_{L^{2j}}\le C\sqrt{2j}\|v\|_{H^1}$.

Define the transport commutator
\[
 \mathfrak C(u,f)=u\cdot\nabla Qf-Q\diver(f\otimes u),
 \qquad (\diver(f\otimes u))_j=\partial_i(f_ju_i).
\]
For each fixed $m>2$,
\begin{equation}\label{pcper:commutator}
 \|\mathfrak C(u,f)\|_{L^m}
                  \le C_m\|\nabla u\|_{L^2}\|f\|_{L^m}.
\end{equation}
This is the periodic commutator estimate of~\cite{CRW1976},
with $\|u\|_{\mathrm{BMO}}\le C\|\nabla u\|_{L^2}$; see also
\cite[proof of Lemma~2.3]{SourceNS}.

\subsection{Density moments}

Put $\psi=Q(\rho u)$ and $\mathcal K=Q(\rho\nabla\Phi)$.
The primitive transport equation is
\begin{equation}\label{pcper:theta-transport}
 D_t(\theta(\rho)-\psi)+P
 =1-\bar G-\mathfrak C(u,\rho u)-\mathcal K.
\end{equation}
Indeed, application of $Q$ to conservative momentum gives
$D_t\psi=-G+\bar G+\mathfrak C(u,\rho u)+\mathcal K$,
whereas $D_t\theta+P=1-G$.

Choose a fixed integer $k$ such that
\begin{equation}\label{pcper:k-choice}
 k\ge7,\qquad k\ge2\gamma+1,\qquad (2\kappa-1)k\ge1,
 \qquad m=k+1.
\end{equation}
For a fixed shift $L>0$, write
\[
 \varphi=(\theta(\rho)-\psi-L)_+,\quad
 \mathcal M_k=\int\rho\varphi^k\dd,\quad X=1+\mathcal M_k,
 \quad Y=\int\rho\theta(\rho)_+^k\dd.
\]
The mass and kinetic energy, \eqref{pcper:potential-estimates}, and
mass interpolation give, with $r=2m/(m+1)$,
\[
 \|\psi\|_{L^{2m}}\le C_m\|\rho u\|_{L^r}
       \le C_m\|\rho\|_{L^m}^{1/2},\qquad
 \|\rho\|_{L^r}\le\|\rho\|_{L^m}^{1/2},\qquad
 \int\rho|\psi|^k\dd\le C_m\|\rho\|_{L^m}^{m/2}.
\]
By \eqref{pc:growth},
$\int\rho^m[\log(e+\rho)]^{\kappa k}\le C(1+Y)$.
Using $\theta_+\le\varphi+|\psi|+L$ and absorbing the resulting
square root proves
\begin{equation}\label{pcper:moment-comparison}
\begin{gathered}
 \int\left(\rho\theta(\rho)_+^k+
              \rho^m[\log(e+\rho)]^{\kappa k}\right)\dd\le C_{k,L}X,\\
 \|\rho\|_{L^m}\le C_{k,L}X^{1/m},\quad
 \|\psi\|_{L^{2m}}\le C_{k,L}X^{1/(2m)},\quad
 \int\rho|\psi|^k\dd\le C_{k,L}X^{1/2},\\
 \|\psi\|_{L^\infty}\le C\|\rho u\|_{L^{8/3}}
       \le C\|\rho\|_{L^4}\|u\|_{L^8}\le C_{k,L}A_3X^{1/m}.
\end{gathered}
\end{equation}

\begin{lemma}[Electric absorption]\label{pcper:electric-absorption-lemma}
For each $\varepsilon>0$, one can choose a fixed $L$ sufficiently
large, independently of the density maximum and of $S$, so that
\begin{equation}\label{pcper:electric-absorption}
 \int\rho\varphi^{k-1}|\mathcal K|\dd
 \le\varepsilon\int\rho^{\gamma+1}\varphi^{k-1}\dd
                         +C_{k,L,\varepsilon}A_3^2X.
\end{equation}
\end{lemma}
\begin{proof}
Since $\gamma+1>2$, fixed-exponent Poisson estimates imply
\begin{equation}\label{pcper:electric-potential-infinity}
 \|\mathcal K\|_{L^\infty}
 \le C_\gamma\|\rho\|_{L^{\gamma+1}}\|\nabla\Phi\|_{L^\infty}
 \le C_\gamma(1+\|\rho\|_{L^{\gamma+1}}^2).
\end{equation}
Put $g(\rho)=(\theta(\rho)-L)_+^{k-1}$,
$W_0=\int\rho g(\rho)$ and $D_0=\int\rho^{\gamma+1}g(\rho)$.
When $W_0>0$, monotonicity of $r^\gamma$ and $g(r)$ in the
probability measure $\mathrm d\nu=\rho\dd$ gives
\[
 \int\rho^\gamma g(\rho)\,\mathrm d\nu
 -\int\rho^\gamma\,\mathrm d\nu\int g(\rho)\,\mathrm d\nu
 =\frac12\iint(\rho(x)^\gamma-\rho(y)^\gamma)
       (g(\rho(x))-g(\rho(y)))\,\mathrm d\nu(x)\,\mathrm d\nu(y)\ge0.
\]
Thus, writing $A=D_0/W_0$,
\[
 \|\rho\|_{L^{\gamma+1}}^2\le A^{2/(\gamma+1)},\qquad
 A\ge[\theta^{-1}(L)]^\gamma,
 \qquad
 \frac{\|\mathcal K\|_{L^\infty} W_0}{D_0}
 \le C_\gamma\bigl(A^{-1}+A^{-(\gamma-1)/(\gamma+1)}\bigr).
\]
Choose $L$ so large that this ratio is at most a prescribed
$\delta>0$. The same inequality is trivial if $W_0=0$. Consequently,
\begin{equation}\label{pcper:unshifted-absorption}
 \|\mathcal K\|_{L^\infty} W_0\le\delta D_0.
\end{equation}

Since $\gamma+1\le m$, \eqref{pcper:moment-comparison} and
\eqref{pcper:electric-potential-infinity} give
$\|\mathcal K\|_{L^\infty}\le C_{k,L}X^{2/m}$ and
$\int\rho|\psi|^{k-3}\le C_{k,L}X^{1/2}$.
It follows that
\begin{equation}\label{pcper:electric-shift-errors}
\begin{aligned}
 \|\mathcal K\|_{L^\infty}\int\rho|\psi|^{k-1}\dd
 &\le C_{k,L}A_3^2X^{1/2+4/m}\le C_{k,L}A_3^2X,\\
 \int\rho^{\gamma+1}|\psi|^{k-1}\dd
 &\le\|\psi\|_{L^\infty}^2\|\rho\|_{L^s}^{\gamma+1}
                         \|\psi\|_{L^{2m}}^{k-3}
 \le C_{k,L}A_3^2X^{(k+2\gamma+3)/(2m)}
 \le C_{k,L}A_3^2X.
\end{aligned}
\end{equation}
Here $s=2m(\gamma+1)/(k+5)$ satisfies $1<s\le m$ and
$(\gamma+1)/s+(k-3)/(2m)=1$.
The last two absorptions use respectively $m\ge8$ and
$k\ge2\gamma+1$.
Finally, with $c_k=2^{k-2}$,
\[
 \varphi^{k-1}\le c_k(g(\rho)+|\psi|^{k-1}),\qquad
 g(\rho)\le c_k(\varphi^{k-1}+|\psi|^{k-1}).
\]
Apply \eqref{pcper:unshifted-absorption} with
$\delta=\varepsilon/c_k^2$ and then
\eqref{pcper:electric-shift-errors}. This proves
\eqref{pcper:electric-absorption}.
\end{proof}

\begin{lemma}[Uniform moments and tails]\label{pcper:moments-lemma}
For each $k$ satisfying \eqref{pcper:k-choice},
\begin{equation}\label{pcper:moments}
 \sup_{t\le S}\int\left(\rho\theta(\rho)_+^k+
       \rho^{k+1}[\log(e+\rho)]^{\kappa k}\right)\dd\le C_k.
\end{equation}
Every fixed finite density moment is uniformly bounded, and
\begin{equation}\label{pcper:momentum-moments}
 \|\rho u\|_{L^p}\le C_p A_3\quad(2\le p<\infty).
\end{equation}
For each fixed $p\ge1$, there are fixed $L_p,C_p$ such that
\begin{equation}\label{pcper:tails}
 \int_0^S\int(\rho-L_p)_+^p\dd\,\mathrm dt\le C_p.
\end{equation}
\end{lemma}
\begin{proof}
The zero-momentum coercivity and \eqref{pcper:GN} imply
$\int\exp(c|u|^2/A_3^2)\le C$ when $A_3>0$;
$u=0$ when $A_3=0$. The Orlicz--Young inequality applied to
$\rho^m$ and $(c_*|u|/A_3)^m$ gives
\[
 \|\rho u\|_{L^m}\le C_{k,L}A_3X^{1/m},
 \qquad \kappa k\ge m/2.
\]
For a scalar $Z$, H\"older in the mass measure gives
\[
 \int\rho\varphi^{k-1}|Z|\dd
 \le\mathcal M_k^{1-1/k}\|\rho\|_{L^m}^{1/k}\|Z\|_{L^m}.
\]
Together with \eqref{pcper:commutator}, the commutator contribution
is at most $C_{k,L}A_3^2X$, because
$1-1/k+1/(km)+1/m=1$.
Since $\bar G=\int\lambda(\rho)\diver u-(\bar P-1)$ and
$\bar P$ is bounded by the energy, it suffices to estimate
$V=\int\lambda(\rho)\diver u$. By \eqref{pc:coefficient},
\[
 |V|^2\le A_3^2\int\lambda(\rho)\dd
       \le C_{k,L}A_3^2X^{1/k}.
\]
Using $|V|\le |V|^2+1$ and
$\int\rho\varphi^{k-1}\le\mathcal M_k^{1-1/k}$, testing
\eqref{pcper:theta-transport} in the transported mass measure and
applying \eqref{pcper:electric-absorption} with a small fixed
constant gives
\[
 \mathcal M_k'+\frac{3k}{4}\int\rho^{\gamma+1}\varphi^{k-1}\dd
 \le C_{k,L}A_3^2X+Ck\int\rho\varphi^{k-1}\dd.
\]
The last constant $C$ is independent of $L$. Increase $L$ so that
$L>1+\max\{\theta((8C)^{1/\gamma}),0\}$. On
$\{\rho^\gamma\ge8C\}$ the pressure absorbs the last integral.
On its complement, $\varphi\le|\psi|$, and
\[
 \int_{\{\rho^\gamma<8C\}}\rho\varphi^{k-1}\dd
 \le\|\psi\|_{L^\infty}^2\int\rho|\psi|^{k-3}\dd
 \le C_{k,L}A_3^2X^{2/m+1/2}\le C_{k,L}A_3^2X.
\]
We have therefore proved
\begin{equation}\label{pcper:moment-differential}
 \mathcal M_k'+c_k\int\rho^{\gamma+1}\varphi^{k-1}\dd
                   \le C_{k,L}A_3^2(1+\mathcal M_k).
\end{equation}
The fixed initial $L^\infty$ and $H^1$ bounds imply
$\mathcal M_k(0)\le C_k$. Gronwall and
\eqref{pcper:coercivity} give
\begin{equation}\label{pcper:moment-damping}
 \sup_{t\le S}\mathcal M_k(t)
 +\int_0^S\int\rho^{\gamma+1}\varphi^{k-1}\dd\,\mathrm dt\le C_k.
\end{equation}
Equation \eqref{pcper:moment-comparison} proves
\eqref{pcper:moments}; $k$ can be arbitrarily large, which gives
all fixed density moments and \eqref{pcper:momentum-moments}.

For $p\ge3$, choose an admissible $k\ge p$ and then $L_p$ so large
that $\theta(r)\ge2L$ and $\theta(r)\ge cr$ for $r\ge L_p$.
On that set either $\varphi\ge\theta(\rho)/4$ or
$\psi\ge\theta(\rho)/4\ge c\rho$. Thus
\[
 \int_{\{\rho\ge L_p\}}\rho^p\dd
 \le C_p\int\rho^{\gamma+1}\varphi^{k-1}\dd
     +C_p\|\psi\|_{L^\infty}^2\int\rho|\psi|^{p-3}\dd
 \le C_p\int\rho^{\gamma+1}\varphi^{k-1}\dd+C_pA_3^2.
\]
Here the established fixed moments give $\|\psi\|_{L^\infty}\le C_pA_3$
and a uniform bound for the last weighted potential moment.
Integration proves \eqref{pcper:tails}. For $1\le p<3$, increase
$L_3$ by one and use $(\rho-L_3-1)_+^p\le(\rho-L_3)_+^3$.
\end{proof}

\subsection{Pressure dissipation and electric work}

\begin{lemma}[Uniform pressure budget]\label{pcper:pressure-lemma}
Define
\begin{equation}\label{pcper:budget-definition}
 D=A_3^2+\|P-1\|_{L^2}^2+\|P-1\|_{L^4}^4.
\end{equation}
Then
\begin{equation}\label{pcper:pressure-budget}
 \int_0^S\left[D+
       \int(\rho+1)^{\gamma-1}(\rho-1)^2\dd
             +\int\rho|\nabla\Phi|^2\dd\right]\mathrm dt\le C.
\end{equation}
\end{lemma}
\begin{proof}
Put $J=\int\rho u\cdot\nabla\Phi$. Testing conservative momentum
by $\nabla\Phi$ gives the exact identity
\begin{equation}\label{pcper:pressure-test}
\begin{aligned}
 \int (P-1)(\rho-1)\dd+\int\rho|\nabla\Phi|^2\dd
 ={}&J'+\|\nabla\Delta^{-1}\diver(\rho u)\|_{L^2}^2
       -\int\rho u\otimes u:\nabla^2\Phi\dd\\
 &+\int h\diver u\,(\rho-1)\dd.
\end{aligned}
\end{equation}
The time derivative of the test contributes
$\int\rho u\cdot\nabla\Phi_t=-\|\nabla\Phi_t\|_{L^2}^2$.
The fixed moments give
\[
 \|\nabla\Phi_t\|_{L^2}^2\le\|\rho u\|_{L^2}^2\le CA_3^2,
 \quad
 \left|\int\rho u\otimes u:\nabla^2\Phi\dd\right|
       \le C\|u\|_{L^8}^2\le CA_3^2,\quad |J|\le C.
\]
By \eqref{pc:coefficient}, split the viscosity term on $\{\rho\le2L\}$ and its complement, with
$L\ge L_{\ell+2}$, $L\ge2$. Weighted Cauchy--Schwarz and Young give
\[
 \left|\int h\diver u\,(\rho-1)\dd\right|
 \le\varepsilon\|\rho-1\|_{L^2}^2+C_\varepsilon A_3^2
                         +C\int(\rho-L)_+^{\ell+2}\dd.
\]
Since $(P-1)(\rho-1)\ge c_\gamma(\rho+1)^{\gamma-1}(\rho-1)^2$
and this dominates $c_\gamma|\rho-1|^2$, integration of
\eqref{pcper:pressure-test} and \eqref{pcper:tails} proves the
weighted density and electric assertions. For $p=2,4$, increasing
a fixed threshold if necessary gives
\[
 |\rho^\gamma-1|^p\le C_p|\rho-1|^2
                         +C_p(\rho-L_{p\gamma})_+^{p\gamma}.
\]
This proves the rest of \eqref{pcper:pressure-budget}.
\end{proof}

The fixed moments and periodic elliptic regularity give
\begin{equation}\label{pcper:electric-work}
\begin{gathered}
 |\rho-1|\le|P-1|,\qquad
 \|\rho\nabla\Phi\|_{L^2}+\|\sqrt\rho\nabla\Phi\|_{L^2}
                                   \le C\|P-1\|_{L^2},\\
 \left|\int\rho\nabla\Phi\cdot\dot u\dd\right|
                          \le\varepsilon A_2^2+C_\varepsilon\|P-1\|_{L^2}^2,\\
 \|D_t\nabla\Phi\|_{L^2}
 \le C\|\rho u\|_{L^2}+\|u\|_{L^4}\|\nabla^2\Phi\|_{L^4}\le CA_3,\\
 \left|\int\rho D_t\nabla\Phi\cdot\dot u\dd\right|
 \le\|\rho\|_{L^4}\|D_t\nabla\Phi\|_{L^2}\|\dot u\|_{L^4}
 \le\varepsilon\|\nabla\dot u\|_{L^2}^2+C_\varepsilon A_3^2.
\end{gathered}
\end{equation}
The first line follows from
$\|\nabla\Phi\|_{L^4}\le C\|\rho-1\|_{L^2}$ and the fixed density
moments; the last uses \eqref{pcper:acceleration-mean}.

\subsection{Flux and acceleration estimates}

Fix the provisional density bound $R=R_S$.
The relative momentum identity and its elliptic consequences are
\begin{equation}\label{pcper:flux-elliptic}
\begin{gathered}
 \rho\dot u=\nabla G+\mu\nabla^\perp\omega+\rho\nabla\Phi,\\
 \|\nabla G\|_{L^2}+\|\nabla\omega\|_{L^2}
              \le C(\sqrt R\,A_2+\|P-1\|_{L^2}),\\
 |\bar G|\le C(A_3+\|P-1\|_{L^2}),\qquad
 \|G\|_{L^2}\le CR^{\ell/2}A_1.
\end{gathered}
\end{equation}
Indeed $\bar G=\int\lambda\diver u-(\bar P-1)$ and
$\int\lambda\le C$ by the fixed moments.
The definitions and the same moments imply
\begin{equation}\label{pcper:energy-comparison}
\begin{gathered}
 A_1^2\le C(A_3^2+\|P-1\|_{L^2}^2)\le CD,\qquad
 A_3^2\le C(A_1^2+\|P-1\|_{L^2}^2),\\
 \mathscr E\asymp1+A_3^2,\qquad D\le C\mathscr E,
 \qquad \sup_{t\le S}(\|P-1\|_{L^2}+\|P-1\|_{L^4})\le C.
\end{gathered}
\end{equation}

\begin{lemma}[Logarithmic estimates]\label{pcper:logarithmic-lemma}
There exist $c>0$ and a fixed finite $N\ge1$ such that
\begin{equation}\label{pcper:log-differential}
 \mathscr E'+cA_2^2\le CR\log(e+R)\mathscr E D,
\end{equation}
and
\begin{equation}\label{pcper:log-energy}
 \sup_{t\le S}\log\mathscr E(t)
       +\int_0^S\frac{A_2^2}{\mathscr E}\,\mathrm dt\le CR\log(e+R).
\end{equation}
Moreover,
\begin{equation}\label{pcper:acceleration-differential}
 (A_2^2)'+c\|\nabla\dot u\|_{L^2}^2\le CR^ND(e+A_2^2),
\end{equation}
and, for $0<\tau<S$,
\begin{equation}\label{pcper:log-acceleration}
 \sup_{\tau\le t\le S}\log(e+A_2^2(t))
       +\int_\tau^S\frac{\|\nabla\dot u\|_{L^2}^2}{e+A_2^2}\,\mathrm dt
                  \le C_\tau R^N.
\end{equation}
\end{lemma}
\begin{proof}
Differentiation gives
\[
 D_tG=h\diver\dot u-\rho\lambda'\,(\diver u)^2
           -h\operatorname{tr}((\nabla u)^2)+\gamma P \diver u,
 \qquad D_t\omega=\curl\dot u-\omega \diver u.
\]
Testing momentum by $2\dot u$, with no boundary terms, yields
\begin{equation}\label{pcper:exact-flux-energy}
\begin{aligned}
 \mathscr E'+2A_2^2={}&
 \int\left(1+\frac{\rho\lambda'}h\right)\frac{G^2}h \diver u\dd
 -2\int\left(1+\frac{\rho\lambda'}h\right)G(\diver u)^2\dd\\
 &+2\gamma\int\frac{GP}h \diver u\dd
 +4\int G\det\nabla u\dd-\mu\int\omega^2\diver u\dd\\
 &
                        +2\int\rho\nabla\Phi\cdot\dot u\dd.
\end{aligned}
\end{equation}
Put $d=A_3+\|P-1\|_{L^2}$. For $0<\delta\le1/4$, weighted H\"older,
\eqref{pcper:GN}, and \eqref{pcper:flux-elliptic} give
\begin{equation}\label{pcper:weighted-quadratic}
 \left\|\frac{G^2}h\right\|_{L^2}
 \le\frac C{\sqrt\delta}R^{\ell\delta/2}
                        A_1(\sqrt R\,A_2+d)
                   +CA_1^{1-\delta}d^{1+\delta}.
\end{equation}
In detail the first step is
$\|G^2/h\|_{L^2}\le CA_1^{1-\delta}
\|G\|_{L^{2(1+\delta)/\delta}}^{1+\delta}$, and
\[
 \|G-\bar G\|_{L^{2(1+\delta)/\delta}}^{1+\delta}
 \le C\delta^{-1/2}\|G-\bar G\|_{L^2}^\delta\|\nabla G\|_{L^2}.
\]
The uniformity follows from
$\delta^{-(1+\delta)/2}\le C\delta^{-1/2}$.
The absolute sum of the first three integrals in
\eqref{pcper:exact-flux-energy} is at most
\[
 CA_3\left(\|G^2/h\|_{L^2}+\|G(P+1)/h\|_{L^2}\right),
 \qquad \|G(P+1)/h\|_{L^2}\le C\|G^2/h\|_{L^2}^{1/2}.
\]
Here the second inequality uses $h\ge2\mu$ and the fixed $L^4$
pressure moment. The $A_2^2$ term in the quadratic expression costs
$\varepsilon A_2^2+C_\varepsilon\delta^{-1}R^{1+\ell\delta}A_3^2A_1^2$.
Its square-root counterpart costs
\[
 \varepsilon A_2^2+C_\varepsilon\delta^{-1/3}R^{(1+\ell\delta)/3}
                                      A_3^{4/3}A_1^{2/3}.
\]
The last energy product is at most $CD$. All remaining products
are bounded by $C(1+\delta^{-1/2}R^{\ell\delta/2})\mathscr E D$,
using $A_1^2+d^2\le CD\le C\mathscr E$.
Choose $\delta=[4\log(e+R)]^{-1}$.

The Hardy--BMO Jacobian estimate~\cite{CLMS1993}, used in
\cite[proof of Lemma~2.6]{SourceNS}, gives
\[
\begin{aligned}
 \left|\int G\det\nabla u\dd\right|
 &\le C\|\nabla G\|_{L^2}\|\nabla u\|_{L^2}^2\\
 &\le C(\sqrt R\,A_2+d)A_3^2
 \le\varepsilon A_2^2+C_\varepsilon R\mathscr E D.
\end{aligned}
\]
The mean of $G$ drops out because $\int\det\nabla u\dd=0$.
The vorticity term satisfies
\[
 \left|\int\omega^2\diver u\dd\right|
 \le C A_3^2(\sqrt R\,A_2+d)
                   \le\varepsilon A_2^2+C_\varepsilon R\mathscr E D.
\]
Finally use \eqref{pcper:electric-work}. This proves
\eqref{pcper:log-differential}; division by $\mathscr E$ and
\eqref{pcper:pressure-budget} prove \eqref{pcper:log-energy}.
The initial upper bound for $A_1$ follows from $C_0$ and the
coefficient bounds on the initial density range.

We next give the material derivative calculation. Let
$J_0=\left(\begin{smallmatrix}0&1\\-1&0\end{smallmatrix}\right)$,
$\mathsf S=(G-1)I+\mu\omega J_0$ and
$\Sigma=D_t\mathsf S+(\diver u)\mathsf S-\mathsf S(\nabla u)^{\mathsf T}$.
The continuity equation and divergence commutation give
\begin{equation}\label{pcper:commuted-stress}
\begin{gathered}
 \rho D_t\dot u=\diver\Sigma+\rho D_t\nabla\Phi,
 \\
 \Sigma=h(\diver\dot u)I+\mu(\curl\dot u)J_0+\mathcal Q,
\end{gathered}
\end{equation}
where
\[
\begin{aligned}
 \mathcal Q={}&\bigl((h-\rho\lambda')(\diver u)^2
       -h\operatorname{tr}((\nabla u)^2)+(\gamma-1)P\diver u\bigr)I\\
       &-(G-1)(\nabla u)^{\mathsf T}-\mu\omega J_0(\nabla u)^{\mathsf T},\\
 |\mathcal Q|\le{}&C(h|\nabla u|^2+P|\nabla u|).
\end{aligned}
\]
Testing by $\dot u$ and applying \eqref{pcper:electric-work} yields
\[
 (A_2^2)'+c\|\nabla\dot u\|_{L^2}^2
 \le CR^{2\ell}\|\nabla u\|_{L^4}^4
                   +C\int P^2|\nabla u|^2\dd+CA_3^2.
\]
To bound these products, first test momentum by $u$:
\[
 A_3^2=\int (P-1)\diver u\dd-\int\rho\dot u\cdot u\dd
                      +\int\rho\nabla\Phi\cdot u\dd
             \le CA_3+CA_2.
\]
Hence $\mathscr E^2\le C(1+A_2^2)$.
Fixed-exponent Hodge and Gagliardo--Nirenberg, with
\eqref{pcper:flux-elliptic}, give
\begin{equation}\label{pcper:fourth-powers}
\begin{aligned}
 \|\nabla u\|_{L^4}^4
 &\le C R^\ell A_1^2(RA_2^2+d^2)+CA_3^4+C\|P-1\|_{L^4}^4+Cd^4\\
 &\le CR^ND(e+A_2^2),\qquad \mathscr E^2\le C(1+A_2^2).
\end{aligned}
\end{equation}
Also
$\int P^2|\nabla u|^2\le C(A_3^2+\|P-1\|_{L^4}^4+
\|\nabla u\|_{L^4}^4)$. This proves
\eqref{pcper:acceleration-differential}, after increasing $N$.
From \eqref{pcper:log-energy},
\[
 \int_0^S A_2^2\,\mathrm dt
 \le CR\log(e+R)\exp(CR\log(e+R)).
\]
Averaging on $[\tau/2,\tau]$ supplies a starting point with
$\log(e+A_2^2)\le C_\tau R^N$. Divide
\eqref{pcper:acceleration-differential} by $e+A_2^2$, integrate,
and use \eqref{pcper:pressure-budget}. No initial acceleration
norm is required. This proves \eqref{pcper:log-acceleration}.
\end{proof}

\subsection{Uniform density bound}

\begin{lemma}[Uniform flux cost]\label{pcper:flux-lemma}
For every $0<\tau<S$,
\begin{equation}\label{pcper:weighted-flux}
 \int_\tau^S\frac{\|G\|_{L^\infty}^2}{\mathscr E}\,\mathrm dt
                      \le C_\tau R_S^2[\log(e+R_S)]^2.
\end{equation}
On every prefix interval $[0,t]\subset[0,S]$ and $\tau\le s<t$,
\begin{equation}\label{pcper:local-flux}
 \int_s^t\|G\|_{L^\infty}\,\mathrm d\xi
 \le C_\tau R_t\log(e+R_t)
            \left(t-s+\int_s^t A_3^2\,\mathrm d\xi\right)^{1/2}.
\end{equation}
\end{lemma}
\begin{proof}
For every $m\ge2$, the zero-mean field satisfies
\[
 \|\nabla\Phi\|_{L^m}\le C\sqrt m\,\|\rho-1\|_{L^2}.
\]
This is the $H^1$ Sobolev estimate applied to $\nabla\Phi$;
it does not require a density moment whose exponent grows with $m$.
Consequently, periodic Calder\'on--Zygmund and
\eqref{pcper:acceleration-mean} give
\begin{equation}\label{pcper:high-gradient}
\begin{aligned}
 \|\nabla G\|_{L^m}
 &\le Cm\bigl(\|\rho\dot u\|_{L^m}
                         +\|\rho\nabla\Phi\|_{L^m}\bigr)\\
 &\le Cm^{3/2}R\bigl(\|\nabla\dot u\|_{L^2}+\|P-1\|_{L^2}\bigr).
\end{aligned}
\end{equation}
Lemma~2.9 in~\cite{SourceNS} gives, for $\bar f=0$ and $p,m\ge8$,
\begin{equation}\label{pcper:log-interpolation}
\begin{gathered}
 \|f\|_{L^\infty}\le C\sqrt p\,
 \|f\|_{L^2}^\eta\|\nabla f\|_{L^2}^{1-\eta-\zeta}
                         \|\nabla f\|_{L^m}^\zeta,
 \\
 \eta=\frac{2(m-2)}{p(m-2)+2m},\quad
 \zeta=\frac{2m}{p(m-2)+2m}.
\end{gathered}
\end{equation}
The constant is independent of $p,m$.

Put
\[
 a_0=\frac{A_1}{\sqrt{\mathscr E}},\quad
 L_0=\frac{A_2+A_3+\|P-1\|_{L^2}}{\sqrt{\mathscr E}},\quad
 Z=\frac{\|\nabla\dot u\|_{L^2}}{\sqrt{e+A_2^2}}.
\]
The already proved estimates imply, for a fixed sufficiently large $N$,
\begin{equation}\label{pcper:flux-time-factors}
\begin{gathered}
 a_0\le1,\qquad \int_0^S a_0^2\,\mathrm dt\le C,
 \qquad \int_0^S L_0^2\,\mathrm dt\le CR\log(e+R),\\
 \int_\tau^S Z^2\,\mathrm dt\le C_\tau R^N,
 \qquad \sup_{[\tau,S]}\log(e+L_0^2)\le C_\tau R^N.
\end{gathered}
\end{equation}
Choose $p=8+(N+\ell+4)\log(e+R)$ and $m=e^p$.
Apply \eqref{pcper:log-interpolation} to $G-\bar G$.
Equations \eqref{pcper:flux-elliptic} and
\eqref{pcper:high-gradient} imply
\[
\begin{aligned}
 \frac{\|G-\bar G\|_{L^\infty}}{\sqrt{\mathscr E}}
 &\le C\sqrt p\,
 R^{1/2+(\ell-1)\eta/2+\zeta/2}m^{3\zeta/2}
 a_0^\eta L_0^{1-\eta-\zeta}
       (L_0+\mathscr E^{-1/2}Z+L_0Z)^\zeta\\
 &\le C\sqrt{pR}\,a_0^\eta L_0^{1-\eta-\zeta}
                         (L_0^\zeta+Z^\zeta+L_0^\zeta Z^\zeta).
\end{aligned}
\]
All extra powers of $R$ and $m$ are bounded. H\"older's inequality
and \eqref{pcper:flux-time-factors} bound the $L^2(\tau,S)$ norms
of the first two summands, after removing $C\sqrt{pR}$, by
$C_\tau[R\log(e+R)]^{1/2}R^{N\zeta/2}$.
For the third summand use
\[
 \zeta-\eta=\frac4{p(m-2)+2m},\qquad
 \sup_{[\tau,S]}L_0^{\zeta-\eta}
       \le\exp\left(\frac{C_\tau R^N}{pm}\right)\le C_\tau.
\]
Discard $a_0^\eta\le1$ and apply H\"older with weights
$(1-\zeta,\zeta)$. Its norm is bounded by
\[
 C_\tau\|L_0\|_{L^2}^{1-\zeta}\|Z\|_{L^2}^\zeta
 \le C_\tau[R\log(e+R)]^{1/2}R^{N\zeta/2}.
\]
Since $R^{N\zeta}$ is bounded,
\[
 \left(\int_\tau^S
       \frac{\|G-\bar G\|_{L^\infty}^2}{\mathscr E}\,\mathrm dt\right)^{1/2}
       \le C_\tau\sqrt{pR}\sqrt{R\log(e+R)}
       \le C_\tau R\log(e+R).
\]
The mean contributes at most
$C\int_0^S(A_3^2+\|P-1\|_{L^2}^2)\,\mathrm dt$.
This proves \eqref{pcper:weighted-flux}. Apply it on each prefix
$[0,t]$ and use $\mathscr E\le C(1+A_3^2)$ and Cauchy--Schwarz
on $[s,t]$ to obtain \eqref{pcper:local-flux}.
\end{proof}

\begin{lemma}[Initial density bound]\label{pcper:initial-lemma}
There are $\tau_0\in(0,1]$ and $M_0$, depending only on the fixed
data, such that
\begin{equation}\label{pcper:initial-bound}
 \sup_{0\le t\le\min\{S,\tau_0\}}\|\rho(t)\|_{L^\infty}\le M_0.
\end{equation}
\end{lemma}
\begin{proof}
Take $M_0=2(1+C_0)$ and work before first contact with $M_0$,
with $t\le1$. Equations \eqref{pcper:log-energy} and
\eqref{pcper:acceleration-differential} imply
\[
 \sup_{s\le t}\mathscr E(s)+\int_0^t A_2^2\,\mathrm ds\le C_*,
 \qquad \sup_{0<s\le t}sA_2^2(s)
          +\int_0^t s\|\nabla\dot u\|_{L^2}^2\,\mathrm ds\le C_*.
\]
For the second assertion multiply the differential inequality by
time. Its source $\int_0^t(e+A_2^2)$ is controlled by the first
assertion; the initial weighted term vanishes for a smooth solution.
The high-gradient estimate at $m=4$, the uniform fixed pressure
moment, and Gagliardo--Nirenberg give
\[
\begin{aligned}
 \|G\|_{L^\infty}&\le C_*(1+\|\nabla\dot u\|_{L^2}^{2/3}),\\
 \int_0^t\|G\|_{L^\infty}\,\mathrm ds
 &\le C_*t+C_*\left(\int_0^t s\|\nabla\dot u\|_{L^2}^2\,\mathrm ds\right)^{1/3}
                   \left(\int_0^t s^{-1/2}\,\mathrm ds\right)^{2/3}\\
 &\le C_*t^{1/3}.
\end{aligned}
\]
Along a trajectory, $D_t\theta+P=1-G$, so
$\theta(\rho(t))\le\theta(C_0)+t+C_*t^{1/3}$.
Choose $\tau_0$ so that the right-hand side at $\tau_0$ is less
than $\theta(M_0)$. This excludes first contact.
\end{proof}

\begin{proposition}[Uniform estimates]
\label{pcper:uniform-bound}
For the smooth positive solutions above, with
$\sigma(t)=\min\{1,t\}$,
\begin{equation}\label{pcper:uniform-estimate}
\begin{aligned}
 &\sup_{0\le t\le S}\bigl(\|\rho(t)\|_{L^\infty}
       +\|u(t)\|_{H^1}^2+A_1^2(t)+A_3^2(t)\bigr)
       +\int_0^S A_2^2\,\mathrm dt\\
 &\qquad+\sup_{0<t\le S}\sigma(t)A_2^2(t)
       +\int_0^S\sigma(t)\|\nabla\dot u\|_{L^2}^2\,\mathrm dt\le C.
\end{aligned}
\end{equation}
The velocity in this display is in the zero-momentum frame.
Replacing it by the original velocity changes only the constant.
\end{proposition}
\begin{proof}
We apply the abstract density lemma~\cite[Lemma~2.12]{SourceNS}.
The periodic flow is complete and density maxima are attained.
The continuity equation gives
\[
 D_t\theta(\rho)+\rho^\gamma=1-G\le1+|G|.
\]
The energy bound \eqref{pcper:coercivity}, the initial density bound
\eqref{pcper:initial-bound}, and the local flux estimate
\eqref{pcper:local-flux} verify its remaining hypotheses, with
constants independent of $S$. Hence
\begin{equation}\label{pcper:density-bound}
 \sup_{0\le t\le S}\|\rho(t)\|_{L^\infty}\le C.
\end{equation}
With $R$ fixed, \eqref{pcper:log-energy} and
\eqref{pcper:energy-comparison} imply
\[
 \sup_{t\le S}(A_1^2+A_3^2+\|u\|_{H^1}^2)
                  +\int_0^S A_2^2\,\mathrm dt\le C.
\]
Multiply \eqref{pcper:acceleration-differential} by $\sigma$.
The only initial-layer source is controlled by
$\int_0^{\min\{1,S\}}(e+A_2^2)\le C$. Gronwall and
\eqref{pcper:pressure-budget} give
$\sup\sigma A_2^2+\int_0^S\sigma\|\nabla\dot u\|_{L^2}^2\le C$.
This proves \eqref{pcper:uniform-estimate}.
\end{proof}

\section{Periodic strong solutions and relaxation}\label{pc:periodic-strong}
\subsection{Strong estimates}

We prove Theorem~\ref{pcstr:main} using the uniform estimates of
Proposition~\ref{pcper:uniform-bound}. Throughout this section,
$q>2$ is fixed.

\begin{lemma}[Finite-time estimates]\label{pcstr:strong-estimates}
Let a smooth positive periodic solution on $[0,S]$, $S\le T$, satisfy
\eqref{pcper:uniform-estimate}. For every $1<s<q/(q-1)$,
\begin{equation}\label{pcstr:strong-estimate}
\begin{aligned}
 &\sup_{t\le S}\bigl(\|\rho(t)\|_{W^{1,q}}+
                                      \|\rho_t(t)\|_{L^2}\bigr)
             +\int_0^S\|\nabla u\|_{L^\infty}\,\mathrm dt\\
 &+\sup_{0<t\le S}t\|u(t)\|_{H^2}^2
             +\int_0^S\bigl(\|u\|_{H^2}^2+\|u\|_{W^{2,q}}^s
                               +\|\sqrt\rho\,u_t\|_{L^2}^2\bigr)\,\mathrm dt\\
 &+\int_0^S t\bigl(\|u\|_{W^{2,q}}^2+
                                      \|u_t\|_{H^1}^2\bigr)\,\mathrm dt
       \le C_{q,s,T}.
\end{aligned}
\end{equation}
The constant may depend on $\|\rho_0\|_{W^{1,q}}$ and local
$C^1$ bounds for $\lambda$, but is independent of $S$ and of a
positive density lower bound.
\end{lemma}

\begin{proof}
Work in the zero-momentum frame. The density upper
bound in \eqref{pcper:uniform-estimate} is fixed. Periodic Poisson
regularity gives
\begin{equation}\label{pcstr:electric-regularity}
 \|\nabla\Phi\|_{W^{1,p}}\le C_p\|\rho-1\|_{L^p},\quad 1<p<\infty,
 \qquad \|\rho\nabla\Phi\|_{L^q}\le C_q.
\end{equation}
The bounded density, fixed positive mass, and weighted mean
$\int\rho\dot u=0$ from \eqref{pcper:acceleration-mean} give the
same weighted Poincar\'e estimate as in the proof of
\cite[Lemma~4.1]{SourceNS}. Its time-weighted interpolation argument
uses only these facts, $\int_0^S A_2^2\,\mathrm dt\le C$, and
$\int_0^S t\|\nabla\dot u\|_{L^2}^2\,\mathrm dt\le C_T$,
both supplied by \eqref{pcper:uniform-estimate}. It therefore yields
\begin{equation}\label{pcstr:force-integrability}
 \int_0^S\left(\|\rho\dot u\|_{L^q}^s
              +t\|\rho\dot u\|_{L^q}^2
              +t\|\dot u\|_{H^1}^2\right)\,\mathrm dt\le C_{q,s,T}.
\end{equation}

The effective-flux equation now reads
\begin{equation}\label{pcstr:flux-gradient}
 \nabla G+\mu\nabla^\perp\omega
       =\rho\dot u-\rho\nabla\Phi,\qquad
 \|\nabla G\|_{L^q}+\|\nabla\omega\|_{L^q}
       \le C_q\bigl(\|\rho\dot u\|_{L^q}+1\bigr).
\end{equation}
Set $\Xi=h(\rho)\nabla\rho$. Differentiated continuity gives
\begin{equation}\label{pcstr:weighted-density-gradient}
 \Dt\Xi+(\nabla u)^{\mathsf T}\Xi+\Xi\diver u+\rho\nabla G
                       +\frac{\gamma\rho^\gamma}{h(\rho)}\Xi=0.
\end{equation}
No derivative of $\lambda'$ is present, and the last coefficient
is nonnegative. The gradient-closure argument in the proof of
\cite[Lemma~4.2]{SourceNS} uses this transport identity and the
elliptic flux estimates. Here \eqref{pcstr:electric-regularity}
bounds the additional force $\rho\nabla\Phi$, so
\eqref{pcstr:flux-gradient} gives the same estimates with the force
norm replaced by $1+\|\rho\dot u\|_{L^q}$. The uniform $L^2$ flux
bounds from \eqref{pcper:uniform-estimate} and
\eqref{pcstr:force-integrability} thus give the same
logarithmic Gronwall closure: $\|\nabla\rho\|_{L^q}$ is bounded,
$\|\nabla u\|_{L^\infty}$ is time integrable, and both $W^{2,q}$
time integrals in \eqref{pcstr:strong-estimate} are finite.
The coefficient-gradient absorption at exponent two in that proof
also gives
\begin{equation}\label{pcstr:H2-velocity}
 \|u\|_{H^2}\le C_{q,T}(1+A_2).
\end{equation}
Finally, the convective estimates in the proof of
\cite[Proposition~4.3]{SourceNS} use only the bounded $H^1$ norm
of $u$ and the preceding Hessian bounds. They apply unchanged to
$u_t=\dot u-u\cdot\nabla u$ and
$\rho_t=-u\cdot\nabla\rho-\rho\diver u$, and give the remaining
estimates in \eqref{pcstr:strong-estimate}.
Poisson regularity also gives
\begin{equation}\label{pcstr:electric-strong}
 \|\nabla\Phi\|_{W^{2,q}}\le C_q\|\rho-1\|_{W^{1,q}},\qquad
 \|\nabla\Phi_t\|_{H^1}\le C\|\rho_t\|_{L^2}.
\end{equation}
\end{proof}

\subsection{Existence and uniqueness}

\begin{lemma}[Smooth coefficient approximations]\label{pcstr:coefficient-approximation}
There are smooth nonnegative coefficients
$\lambda_\varepsilon\to\lambda$ locally in $C^1([0,\infty))$
which satisfy \eqref{pc:lambda}--\eqref{pc:growth} with common
constants and a common large-density threshold.
\end{lemma}
\begin{proof}
This is \cite[Lemma~4.4]{SourceNS}, whose coefficient assumptions
are \eqref{pc:lambda}--\eqref{pc:growth}.
\end{proof}

\begin{lemma}[Positive solutions]\label{pcstr:positive-construction}
For smooth coefficients satisfying the preceding structural conditions,
positive $H^2$ density data of mass one and $H^2$ velocity data
generate a global positive strong solution. It satisfies
\eqref{pcper:uniform-estimate} and \eqref{pcstr:strong-estimate}.
\end{lemma}

\begin{proof}
For the local construction, take input velocities with $v(0)=u_0$ in
\[
 v\in C H^2\cap L^2H^3,\qquad
 v_t\in L^\infty L^2\cap L^2H^1
\]
on a short interval $[0,t_*]$. The ball radius is fixed by the positive initial data and coefficient
bounds. Constants in this local construction may depend on the
positive density lower bound; the estimates
\eqref{pcper:uniform-estimate} and \eqref{pcstr:strong-estimate}
used in the vacuum limit are independent of it. Transport $\rho_0$
by $v$, solve the normalized periodic Poisson equation, and solve
\[
 \rho(u_t+v\cdot\nabla u)-\mu\Delta u
       -\nabla((\mu+\lambda(\rho))\diver u)
                 =-\nabla P+\rho\nabla\Phi.
\]
We use the transport--parabolic scheme of
\cite[Section~4.3]{SourceNS}. On a positive coefficient set with
$\|\rho\|_{H^2}\le L$ and
$\|v\|_{L^\infty}\le V$,
\begin{equation}\label{pcstr:local-electric}
\begin{aligned}
 &\|\rho\nabla\Phi\|_{H^1}\le C_L,\qquad
   \|\partial_t(\rho\nabla\Phi)\|_{L^2}
                                      \le C_L\|\rho_t\|_{L^2},\\
 &\|\rho_1\nabla\Phi_1-\rho_2\nabla\Phi_2\|_{L^2}
                               \le C_L\|\rho_1-\rho_2\|_{L^2}.
\end{aligned}
\end{equation}
The first line follows from Poisson regularity and
$\Delta\Phi_t=\rho_t$; all transported densities have mass one.
For the second line write the difference as
$(\rho_1-\rho_2)\nabla\Phi_2+
\rho_1\nabla\Delta^{-1}(\rho_1-\rho_2)$.

The linear estimates are
\[
\begin{aligned}
 \|u\|_{H^2}&\le C_{L,V}(1+\|u_t\|_{L^2}+\|u\|_{L^2}),\\
 \|u\|_{H^3}&\le C_L\bigl[\|u_t\|_{H^1}
                     +(1+\|v\|_{H^2})\|u\|_{H^2}+1\bigr].
\end{aligned}
\]
With $Z=\int\rho(|u|^2+|u_t|^2)\dd$, the energy estimate is
\[
\begin{aligned}
 &Z'+c(\|\nabla u\|_{L^2}^2+\|\nabla u_t\|_{L^2}^2)\\
 &\qquad\le C_{L,V}(1+\|\rho_t\|_{H^1}^2+
                                      \|v_t\|_{L^2}^{4/3})(1+Z).
\end{aligned}
\]
These follow by differentiating the linear equation in time, testing
by $u_t$, and using the transport equation in the mass term.
The differentiated electric force is controlled by
\eqref{pcstr:local-electric}; the viscosity term is bounded using
$\|\lambda'(\rho)\rho_t\diver u\|_{L^2}
\le C_L\|\rho_t\|_{H^1}\|u\|_{H^2}$.
The term containing $v_t$ uses $L^4$ interpolation and Young's
inequality with exponent $4/3$. The spatial estimates follow by
Hodge decomposition and the $H^2$ algebra property in two dimensions.
One can first solve implicit Euler variational problems and apply
these estimates to their time differences. The initial difference
is determined by
\[
 u_t(0)=\rho_0^{-1}
 \bigl[\mu\Delta u_0+\nabla((\mu+\lambda(\rho_0))\diver u_0)
        -\rho_0v(0)\cdot\nabla u_0-\nabla P_0+\rho_0\nabla\Phi_0\bigr]
       \in L^2.
\]
Thus the linear solution has
$u\in C H^2\cap L^2H^3$ and
$u_t\in L^\infty L^2\cap L^2H^1$.

In a ball of radius $K$ in these velocity norms,
\[
 \|\rho(t)\|_{H^2}\le\|\rho_0\|_{H^2}e^{C\sqrt{t_*}K},\qquad
 \|\rho_t\|_{H^1}\le C_LK.
\]
Fix $L=2(1+\|\rho_0\|_{H^2})$ and
$V=1+\|u_0\|_{L^\infty}$. For $K\ge1+\|u_0\|_{H^2}$,
interpolation and $v(0)=u_0$ give
\[
 \|v(t)-u_0\|_{L^\infty}
 \le C\|v(t)-u_0\|_{L^2}^{1/2}\|v(t)-u_0\|_{H^2}^{1/2}
 \le CK\sqrt{t_*}.
\]
Decreasing $t_*$ therefore enforces these fixed $L,V$ bounds and
$\rho\ge\frac12\min\rho_0$. The initial value of $u_t$ is fixed by
the data. Choose $K$ large in terms of these data, and then $t_*$
small; the energy estimate preserves the ball. To check contraction, the transported density
difference satisfies
\[
 \sup_{t\le t_*}\|\rho_1-\rho_2\|_{L^2}^2
       \le C_Kt_*\int_0^{t_*}\|v_1-v_2\|_{H^1}^2\,\mathrm dt.
\]
Test the output velocity difference in $L^2$ and use
\eqref{pcstr:local-electric}. With
$d(w)^2=\sup_{t\le t_*}\|w(t)\|_{L^2}^2+
\int_0^{t_*}\|w\|_{H^1}^2\,\mathrm dt$, this gives
\[
 d(u_1-u_2)^2\le C_K(\sqrt{t_*}+t_*)d(v_1-v_2)^2.
\]
The electric force contributes at most
$\eps\|u_1-u_2\|_{H^1}^2+
C_\eps\|\rho_1-\rho_2\|_{L^2}^2$.
Weak compactness makes the ball complete in this weaker metric,
so iteration yields the local positive solution.

Suppose its maximal $H^2$ time $T^*$ were finite. Apply
\eqref{pcper:uniform-estimate} and Lemma~\ref{pcstr:strong-estimates}
with $q=4$ uniformly on $[0,S]$, $S<T^*$. The flow formula gives
\[
 \inf_{t<T^*,x}\rho(t,x)
 \ge\min\rho_0\exp\left(-\int_0^{T^*}
                           \|\diver u\|_{L^\infty}\,\mathrm dt\right)>0.
\]
Constants $C_*$ in the remainder of this continuation step may depend
on this lower bound and the higher initial norms. Since
$A_2(0)<\infty$, \eqref{pcper:acceleration-differential} and
\eqref{pcper:pressure-budget} imply
\[
 \sup_{t<S}A_2^2+
         \int_0^S\|\dot u\|_{H^1}^2\,\mathrm dt\le C_*.
\]
The additional elliptic source obeys
\[
 \|\rho\nabla\Phi\|_{H^1}
       \le C_*(1+\|\nabla\rho\|_{L^2})\le C_*.
\]
Thus \eqref{pcstr:flux-gradient}, its spatial derivative, and
\eqref{pcstr:H2-velocity} give
\[
 \sup_{t<S}\|u\|_{H^2}+
 \int_0^S\bigl(\|\nabla^2G\|_{L^2}^2+
                              \|\nabla^2\omega\|_{L^2}^2\bigr)\,\mathrm dt
       \le C_*.
\]
Differentiate \eqref{pcstr:weighted-density-gradient} and use
$\nabla\rho\in L^\infty L^4$. The needed product estimates are
\[
 \|\nabla\rho\nabla G\|_{L^2}
       \le C_*(1+\|\nabla^2G\|_{L^2}^{1/2}),\qquad
 \|\nabla\rho\nabla^2u\|_{L^2}
       \le C_*(1+\|u\|_{H^3}^{1/2}).
\]
Since
$\nabla\Xi=h(\rho)\nabla^2\rho+
\lambda'(\rho)\nabla\rho\otimes\nabla\rho$,
the differentiated divergence relation and these products yield
\[
\begin{aligned}
 \|\rho\|_{H^2}&\le C_*(1+\|\nabla\Xi\|_{L^2}),\\
 (\|\nabla\Xi\|_{L^2})'
 &\le C_*(1+\|\nabla u\|_{L^\infty})(1+\|\nabla\Xi\|_{L^2})
           +C_*(\|\nabla^2G\|_{L^2}+\|\nabla^2\omega\|_{L^2}).
\end{aligned}
\]
Gronwall bounds $\rho$ in $L^\infty H^2$.
Differentiated Hodge estimates then give $u\in L^2H^3$ and
$u_t=\dot u-u\cdot\nabla u\in L^2H^1$ up to $T^*$.
The twice-differentiated transport equation has integrable
Lipschitz coefficients and $L^1L^2$ sources, so $\rho$ has a strong
$H^2$ terminal trace. The Hilbert triple gives the strong $H^2$
velocity trace. These positive terminal data restart the local
construction, contradicting maximality.
\end{proof}

\begin{proof}[Proof of Theorem~\ref{pcstr:main}]
Let $j_\varepsilon$ be periodic mollifiers and set
\[
 \rho_{0,\varepsilon}=
       \frac{j_\varepsilon*\rho_0+\varepsilon}{1+\varepsilon},\qquad
 u_{0,\varepsilon}=j_\varepsilon*u_0.
\]
These are positive unit-mass data converging in
$W^{1,q}\times H^1$. Use Lemmas~\ref{pcstr:coefficient-approximation}
and~\ref{pcstr:positive-construction} to obtain global positive
approximations. Their low-order bounds are independent of
$\varepsilon$ and of the time interval. Lemma~\ref{pcstr:strong-estimates}
is uniform on every finite interval because the local
$C^1$ bounds for $\lambda_\varepsilon$ are uniform.

Compactness gives, after extraction,
\[
 \rho_\varepsilon\to\rho\quad\hbox{in }C([0,T];C^\alpha(\TT)),
 \quad 0<\alpha<1-2/q,\qquad
 u_\varepsilon\to u\quad\hbox{in }L^2(0,T;H^1).
\]
For velocity compactness first use $[\tau,T]$ and then let
$\tau\downarrow0$, using the common $L^\infty H^1$ bound.
The conservative equation bounds $(\rho_\varepsilon u_\varepsilon)_t$
in $L^\infty H^{-1}$, while the momenta are bounded in
$L^\infty H^1$. The electric force is bounded in $L^\infty L^2$.
Consequently
\[
 \rho_\varepsilon u_\varepsilon\to\rho u
                   \quad\hbox{in }C([0,T];L^2).
\]
The normalized Poisson operator additionally gives
\[
 \nabla\Phi_\varepsilon\to\nabla\Phi
           \quad\hbox{in }C([0,T];W^{1,p}),\quad 1<p<\infty.
\]
Thus $\rho_\varepsilon\nabla\Phi_\varepsilon\to\rho\nabla\Phi$
strongly in $C([0,T];L^2)$, and the limiting electric term is
identified. Uniform convergence of
$\lambda_\varepsilon(\rho_\varepsilon)$ and the strong velocity
convergence identify the stress and convection.

The differentiated continuity equation, with
$\nabla u\in L^1L^\infty$ and
$\rho\nabla\diver u\in L^1L^q$, gives
$\rho\in C([0,T];W^{1,q})$.
For $\delta>0$, $\rho u/\sqrt{\rho+\delta}\in C L^2$ and
\[
 \left\|\frac{\rho u}{\sqrt{\rho+\delta}}-
                       \sqrt\rho\,u\right\|_{L^2}
       \le\sqrt\delta\|u\|_{L^2}.
\]
This proves the weighted velocity trace. Weak compactness on
$[\tau,T]$ identifies $u_t$ and $\dot u$; lower semicontinuity
followed by $\tau\downarrow0$ preserves all time weights.
Equation~\eqref{pcstr:electric-strong} gives the stated potential
regularity and its time continuity. A diagonal extraction over
integer time intervals and $s\uparrow q/(q-1)$ gives one global
solution with \eqref{pcstr:strong-class}.

For uniqueness let $r=\rho_1-\rho_2$, $w=u_1-u_2$, and
$Y=\int\rho_1|w|^2\dd$. The mass-weighted Poincar\'e inequality
and the density difference equation give
\[
 \frac{\mathrm d}{\mathrm dt}\|r\|_{L^2}
 \le C\|\nabla u_2\|_{L^\infty}\|r\|_{L^2}+C\|w\|_{H^1}.
\]
Integrating this norm inequality from $r(0)=0$, and then applying
Cauchy--Schwarz, yields
\begin{equation}\label{pcstr:density-difference}
 \|w\|_{H^1}^2\le C(Y+\|\nabla w\|_{L^2}^2),\qquad
 \|r(t)\|_{L^2}^2
       \le Ct\int_0^t(Y+\|\nabla w\|_{L^2}^2)\,\mathrm ds.
\end{equation}
The electric force difference is
$r\nabla\Phi_2+\rho_1\nabla\Delta^{-1}r$.
It satisfies, since $\int r=0$,
\[
 \left|\int(r\nabla\Phi_2+\rho_1\nabla\Delta^{-1}r)\cdot w\dd\right|
 \le\eps\|\nabla w\|_{L^2}^2+C_\eps(\|r\|_{L^2}^2+Y).
\]
The remaining terms satisfy the Navier--Stokes difference estimate,
using local Lipschitz bounds for $P$ and $\lambda$.
Testing by $w$ yields
\begin{equation}\label{pcstr:velocity-difference}
 Y'+c\|\nabla w\|_{L^2}^2
 \le C(1+\|\nabla u_2\|_{L^\infty})Y
   +C(1+\|\nabla u_2\|_{L^\infty}^2+
                            \|\dot u_2\|_{L^q}^2)\|r\|_{L^2}^2.
\end{equation}
The strong estimates give
\[
 \int_0^T\left[1+\|\nabla u_2\|_{L^\infty}
       +t(1+\|\nabla u_2\|_{L^\infty}^2+
                            \|\dot u_2\|_{L^q}^2)\right]\,\mathrm dt<\infty.
\]
The common initial density and weighted momentum traces imply
$Y(0+)=0$, because
\[
 \|(\sqrt{\rho_1}-\sqrt{\rho_2})u_2\|_{L^2}^2
       \le\|\rho_1-\rho_2\|_{L^\infty}\|u_2\|_{L^2}^2\to0.
\]
Substitute \eqref{pcstr:density-difference} into
\eqref{pcstr:velocity-difference}, integrate from a positive time,
and let it tend to zero. Gronwall applied to
$\sup_{s\le t}Y(s)+\int_0^t\|\nabla w\|_{L^2}^2\,\mathrm ds$
gives $r=w=0$. The normalized potentials then coincide.
\end{proof}

\subsection{Large-time behavior}

\begin{proof}[Proof of \eqref{pcstr:decay}]
Work in the zero-momentum frame. The density upper bound and
the energy identity \eqref{pcper:energy} give
\begin{equation}\label{pcstr:relative-energy}
 \mathcal H'+A_3^2=0,\qquad
 \mathcal H\asymp\|\sqrt\rho\,u\|_{L^2}^2+\|\rho-1\|_{L^2}^2
                         \le C(A_3^2+\|\rho-1\|_{L^2}^2).
\end{equation}
Here $\|\nabla\Phi\|_{L^2}\le C\|\rho-1\|_{L^2}$ and the
relative internal energy is quadratic on the bounded density interval.
Set $K=\int\rho u\cdot\nabla\Phi\dd$. Continuity and momentum give
\begin{equation}\label{pcstr:pressure-test}
\begin{aligned}
 &\int(P-1)(\rho-1)\dd+\int\rho|\nabla\Phi|^2\dd\\
 &=K'+\|\nabla\Delta^{-1}\diver(\rho u)\|_{L^2}^2
       -\int\rho u\otimes u:D^2\Phi\dd
       +\int h(\rho)\diver u(\rho-1)\dd.
\end{aligned}
\end{equation}
The projection and convection terms are bounded by $CA_3^2$.
The last term is at most
$\eps\|\rho-1\|_{L^2}^2+C_\eps A_3^2$.
Moreover $(r^\gamma-1)(r-1)\ge(r-1)^2$ for $r\ge0$ and
$|K|\le C\mathcal H$. Thus for a sufficiently small fixed
$\delta>0$,
\[
 \mathcal H-\delta K\asymp\mathcal H,\qquad
 (\mathcal H-\delta K)'+c(A_3^2+\|\rho-1\|_{L^2}^2)\le0.
\]
Consequently
\begin{equation}\label{pcstr:energy-tail}
 \mathcal H(t)+\|\rho(t)-1\|_{L^2}^2+
                        \int_t^\infty A_3^2\,\mathrm ds\le Ce^{-ct}.
\end{equation}
In particular, bounded density and interpolation give exponential
decay of $\rho-1$ in every finite $L^p$ space, and
\[
 \|P-1\|_{L^2}^2+\|P-1\|_{L^4}^4
                                  \le C\|\rho-1\|_{L^2}^2.
\]

By \eqref{pcper:uniform-estimate}, the quantities $R$, $\mathscr E$,
and $A_2$ are bounded on $[1,\infty)$. Boundedness of the density gives
$D\le C(A_3^2+\|\rho-1\|_{L^2}^2)$. Hence
\eqref{pcper:log-differential} and
\eqref{pcper:acceleration-differential} give
\[
 (A_1^2)'+cA_2^2\le C(A_3^2+\|\rho-1\|_{L^2}^2),
\]
and, for almost every $t\ge1$,
\[
 (A_2^2)'+c\|\nabla\dot u\|_{L^2}^2
                         \le C(A_3^2+\|\rho-1\|_{L^2}^2).
\]
The energy comparison gives
$\int_t^\infty A_1^2\,\mathrm ds\le Ce^{-ct}$.
Average the first differential inequality over starting times in
$[t-1,t]$ and then integrate to infinity. Repeat with the second
inequality and the resulting $A_2^2$ tail. This yields
\begin{equation}\label{pcstr:derivative-decay}
 A_1^2(t)+A_3^2(t)+A_2^2(t)
       +\int_t^\infty\|\nabla\dot u\|_{L^2}^2\,\mathrm ds
       \le Ce^{-ct},\qquad\text{for a.e. }t\ge2.
\end{equation}
The interval $[1,2]$ is covered by the uniform weighted estimates.
For $2<p<\infty$, Hodge and interpolation now give
\[
 \|u\|_{W^{1,p}}
 \le C_p\bigl[A_1^{2/p}(A_2+\|\rho-1\|_{L^2})^{1-2/p}
                  +A_1+A_3+\|P-1\|_{L^p}\bigr].
\]
For $1\le p\le2$ use $\|u\|_{H^1}\le CA_3$.
Poisson regularity gives the asserted field decay; exponents
$p\le2$ follow by inclusion from the $p=2$ estimate.
Since $u_t\in L^2(\tau,T;H^1)$ for every $\tau>0$,
$u\in C([\tau,T];H^1)$. Weak lower semicontinuity in $W^{1,p}$
extends the velocity estimate to every positive time. Density
continuity and the Poisson operator do the same for density and field.
Returning to the original frame proves \eqref{pcstr:decay}.
\end{proof}

\begin{corollary}[Vacuum and growing density gradients]\label{pcstr:vacuum}
If $\rho_0(x_0)=0$ at some point of $\TT$, then for every finite time
\[
 \min_{\TT}\rho(t)=0,\qquad \|\rho(t)-1\|_{L^\infty}\ge1.
\]
For every $2<r\le q$ there are $C_r,c_r>0$ such that
\[
 \|\nabla\rho(t)\|_{L^r}\ge C_r^{-1}e^{c_rt},\qquad t\ge0.
\]
\end{corollary}
\begin{proof}
On finite intervals $\nabla u\in L^1L^\infty$, and the flow formula is
\[
 \rho(t,X(t,a))=\rho_0(a)
       \exp\left(-\int_0^t\diver u(s,X(s,a))\,\mathrm ds\right).
\]
It preserves the vacuum point. Since $\rho-1$ has zero spatial mean,
\[
 1\le\|\rho-1\|_{L^\infty}
 \le C_r\|\nabla\rho\|_{L^r}^{r/[2(r-1)]}
                         \|\rho-1\|_{L^2}^{(r-2)/[2(r-1)]}.
\]
Equation~\eqref{pcstr:energy-tail} proves the claimed exponential
growth, with a change of constants on the bounded initial interval.
\end{proof}

\section{Whole-space a priori estimates}\label{pc:whole-estimates}
\subsection{Relative energy and velocity control}

Write $\delta=\rho-1$. The gradient projection
$\mathbb Q=\nabla\Delta^{-1}\diver$ satisfies
\begin{equation}\label{pcws:field-time}
 E_t=-\mathbb Q(\rho u),\qquad \diver E=\delta,
 \qquad \curl E=0.
\end{equation}
We first derive a priori estimates for smooth positive solutions on
$[0,S]$ in the finite-time strong class, with
$\rho-1\in C([0,S];H^1\cap W^{1,q})$, $q>2$, and
$u\in L^\infty(0,S;H^1)\cap L^1(0,S;W^{1,\infty})$.
Temporary density bounds may be used to justify the calculations.
Unless explicitly indicated by a provisional density maximum,
the estimate constants are independent of these temporary bounds
and of $S$; their initial-data dependence is specified in each result.
The passage to densities allowing vacuum is made in
Section~\ref{pc:whole-strong}.

\begin{lemma}[Relative energy]\label{pcws:energy-lemma}
For $0\le t\le S$,
\begin{equation}\label{pcws:energy}
 \mathcal H(t)+\int_0^t A_3^2\,\mathrm ds=\mathcal H(0),\qquad
 \mathcal H=\int_{\R^2}
 \left(\frac12\rho|u|^2+H(\rho\mid1)+\frac12|E|^2\right)\dd,
\end{equation}
where $A_3^2=\int(h(\rho)(\diver u)^2+\mu\omega^2)\dd$.
In particular,
\begin{equation}\label{pcws:essential-residual}
 \sup_t\left(
 \int_{\{\rho\le2\}}|\rho-1|^2\dd
 +\int_{\{\rho>2\}}\rho^\gamma\dd
 +\|E\|_{L^{2}}^2+\|\sqrt\rho u\|_{L^{2}}^2\right)
 +\int_0^S\|\nabla u\|_{L^{2}}^2\,\mathrm dt\le C\mathcal H(0).
\end{equation}
The constants in this estimate do not use a density upper or lower
bound.
\end{lemma}
\begin{proof}
The transport identity for relative internal energy is
\[
 \frac{\mathrm d}{\mathrm dt}\int H(\rho\mid1)\dd
 =-\int(P-1)\diver u\dd.
\]
Since $\mathbb Q$ is an orthogonal projection and $E$ is a gradient,
\[
 \frac12\frac{\mathrm d}{\mathrm dt}\|E\|_{L^{2}}^2
 =-\int E\cdot\rho u\dd.
\]
These terms cancel those in the kinetic energy identity. The planar
div--curl identity gives the displayed form of $A_3^2$.
Spatial cutoffs, followed by approximation in the indicated energy
spaces, justify the integrations over $\R^2$.
Finally, convexity gives $H(r\mid1)\asymp|r-1|^2$ on $[0,2]$
and $H(r\mid1)\asymp r^\gamma$ on $[2,\infty)$, which proves
\eqref{pcws:essential-residual}.
\end{proof}

\begin{lemma}[Velocity control without a positive density bound]
\label{pcws:weighted-poincare}
If $\int H(\rho\mid1)\dd\le\mathcal H(0)$, then
\begin{equation}\label{pcws:velocity-control}
 \|v\|_{L^{2}}^2\le4\|\sqrt\rho v\|_{L^{2}}^2
       +C\mathcal H(0)\|\nabla v\|_{L^{2}}^2,\qquad v\in H^1(\R^2).
\end{equation}
\end{lemma}
\begin{proof}
On $B=\{\rho<1/2\}$, $H(\rho\mid1)\ge c_\gamma>0$, so
$|B|\le C\mathcal H(0)$. Hence
\[
 \|v\|_{L^{2}}^2
 \le2\|\sqrt\rho v\|_{L^{2}}^2+|B|^{1/2}\|v\|_{L^{4}}^2
 \le2\|\sqrt\rho v\|_{L^{2}}^2
       +C|B|^{1/2}\|v\|_{L^{2}}\|\nabla v\|_{L^{2}}.
\]
Young's inequality proves \eqref{pcws:velocity-control}.
\end{proof}

The kinetic energy term in \eqref{pcws:velocity-control} is retained
in the whole-space estimates.

The Maxwell identity
\begin{equation}\label{pcws:maxwell}
 \rho E=E+\diver\left(E\otimes E-\frac12|E|^2I\right)
\end{equation}
isolates a linear gradient force and a quadratic field stress. The
stress is in $L^1$ by \eqref{pcws:energy}, but the energy identity
alone supplies no time-integrable bound for $\|E\|_{L^{2}}^2$.

\subsection{Uniform density moments and tails}
\label{pcwm:section}

To obtain uniform density moments, we use the screened potential
\[
 \mathcal B=(1-\Delta)^{-1},\qquad Q_s=\mathcal B\diver,
 \qquad \psi=Q_s(\rho u).
\]
The convolution kernel $\mathcal G$ of $\mathcal B$ is positive and
symmetric, has integral one, belongs to every finite $L^r$, and
decays exponentially at infinity. The kernel of $Q_s$ is bounded by
$C|x|^{-1}$ near zero and decays exponentially at infinity.
Consequently, for fixed $1<q<2$ and $r>2$,
\begin{equation}\label{pcwm:bessel-bounds}
 \|Q_sf\|_{L^{2q/(2-q)}}\le C_q\|f\|_{L^{q}},\qquad
 \|Q_sf\|_{L^{\infty}}\le C_r\|f\|_{L^{r}},\qquad
 \|Q_sf\|_{H^1}\le C\|f\|_{L^{2}}.
\end{equation}

\begin{proposition}[Uniform high-density moments and tails]
\label{pcwm:proposition}
Assume \eqref{pc:lambda}--\eqref{pc:growth}. Let a smooth positive
solution exist on $[0,S]$ with the relative energy bound
\eqref{pcws:energy}, and suppose $\|\rho_0\|_{L^{\infty}}$ is bounded.
Choose an integer $k$ such that
\begin{equation}\label{pcwm:k}
 k\ge9,\qquad k\ge2\gamma+5,\qquad (2\kappa-1)k\ge1,
 \qquad m=k+1.
\end{equation}
There is a fixed $L>0$ such that, with
\[
 \varphi=(\theta(\rho)-\psi-L)_+,\qquad
 \mathcal M_k=\int\rho\varphi^k\dd,
\]
one has
\begin{equation}\label{pcwm:uniform-moments}
\begin{aligned}
 &\sup_{0\le t\le S}\int\left[
       \rho\theta(\rho)_+^k
       +\boldsymbol1_{\{\rho>2\}}\rho^{k+1}
                         [\log(e+\rho)]^{\kappa k}\right]\dd\\
 &\qquad+\int_0^S\!\int\rho^{\gamma+1}\varphi^{k-1}
                                      \dd\,\mathrm dt\le C_k.
\end{aligned}
\end{equation}
In particular, for each fixed $2\le p<\infty$,
$\sup_{t\le S}\|\rho(t)-1\|_{L^{p}}\le C_p$. For every fixed
$p\ge1$, there is a threshold $L_p$ such that
\begin{equation}\label{pcwm:tails}
 \int_0^S\!\int(\rho-L_p)_+^p\dd\,\mathrm dt\le C_p.
\end{equation}
The constants depend on the initial relative energy, the initial
density maximum, the fixed viscosity and pressure data, and the
indicated fixed exponents. They are independent of $S$, of the
density maximum on $(0,S]$, and of a positive lower density bound.
\end{proposition}

\begin{proof}
Throughout the proof write
\[
 X=1+\mathcal M_k,\qquad
 Y=\int\rho\theta(\rho)_+^k\dd,\qquad
 D=\int\rho^{\gamma+1}\varphi^{k-1}\dd.
\]
All constants below have the dependence in the proposition. The
proof establishes \eqref{pcwm:moment-ode} and then integrates it
using the energy dissipation. The relative energy gives
\begin{equation}\label{pcwm:high-mass}
 \int_{\{\rho>2\}}\rho\dd\le C,\qquad
 \int_{\{\rho\le2\}}|\rho-1|^2\dd\le C,\qquad
 \|\sqrt\rho u\|_{L^{2}}+\|E\|_{L^{2}}\le C,
 \qquad \int_0^S A_3^2\,\mathrm dt\le C.
\end{equation}
Moreover, \eqref{pcws:velocity-control} and the div--curl identity give
\begin{equation}\label{pcwm:u2}
 \|u\|_{L^{2}}^2\le C(1+A_3^2),\qquad \|\nabla u\|_{L^{2}}\le CA_3.
\end{equation}

\paragraph{Moment comparison.}
The growth of $\theta$ implies
\begin{equation}\label{pcwm:rho-moment-from-y}
 \int_{\{\rho>2\}}\rho^m[\log(e+\rho)]^{\kappa k}\dd
       \le C(1+Y).
\end{equation}
Decompose the momentum as
\[
 \rho u=\sqrt\rho u+(\sqrt\rho-1)\sqrt\rho u.
\]
On $\{\rho\le2\}$, $|\sqrt\rho-1|\le C|\rho-1|$, while on
$\{\rho>2\}$ it is bounded by $C\sqrt\rho$. Thus
\[
 \|\sqrt\rho-1\|_{L^{2m}}\le C(1+Y^{1/(2m)}).
\]
Apply \eqref{pcwm:bessel-bounds} to the first momentum component in
$L^2$ and to the second in $L^{2m/(m+1)}$. This gives
\[
 \|\psi\|_{L^{2m}}\le C(1+Y^{1/(2m)}).
\]
The same bound holds in $L^r$ for $2<r\le2m$. At $r=2$, use
$\sqrt\rho-1\in L^4$, the product in $L^{4/3}$, and the
$L^{4/3}$ convolution bound for the kernel of $Q_s$.

For $2\le j\le k$, split the integral of $\rho|\psi|^j$ at
$\rho=2$. On the high-density part, interpolate between $L^1$ and $L^m$.
By \eqref{pcwm:high-mass},
\[
 \|\rho\boldsymbol1_{\{\rho>2\}}\|_{L^{2m/(2m-j)}}
 \le C\|\rho\boldsymbol1_{\{\rho>2\}}\|_{L^{m}}^{j/(2k)}.
\]
It follows that
\[
 \int\rho|\psi|^j\dd\le C(1+Y^{j/(2k)}).
\]
On a fixed density interval above one, $\theta(\rho)_+^k$ is
bounded by $C_{k,L}|\rho-1|^2$. On the remaining high-density set,
$\theta_+\le\varphi+|\psi|+L$ and the high-density mass is finite.
Consequently,
$Y\le C_{k,L}X+C_{k,L}Y^{1/2}$. Absorbing the square root proves
\begin{equation}\label{pcwm:comparison}
\begin{aligned}
 &Y+\int_{\{\rho>2\}}\rho^m[\log(e+\rho)]^{\kappa k}\dd
                                      \le C_{k,L}X,\\
 &\|\psi\|_{L^{r}}\le C_{k,L,r}X^{1/(2m)},\quad 2\le r\le2m,\qquad
   \int\rho|\psi|^j\dd\le C_{k,L}X^{1/2},\quad 2\le j\le k.
\end{aligned}
\end{equation}
In particular, $\|\rho-1\|_{L^{4}}\le C_{k,L}X^{1/m}$.
For the supremum estimate, write $\rho u=u+(\rho-1)u$.
Since $Q_su=\mathcal B\diver u$, the kernel bounds,
Gagliardo--Nirenberg, and \eqref{pcwm:u2} give
\begin{equation}\label{pcwm:psi-infty}
 \|\psi\|_{L^{\infty}}
 \le CA_3+C\|\rho-1\|_{L^{4}}\|u\|_{L^{8}}
 \le C_{k,L}X^{1/m}(A_3^{3/4}+A_3).
\end{equation}

\paragraph{Commutator estimate.}
For $A_3>0$, the whole-space Trudinger bound is
\begin{equation}\label{pcwm:trudinger}
 \int\left[\exp\left(\frac{c|u|^2}{A_3^2}\right)-1\right]\dd
                  \le C\frac{\|u\|_{L^{2}}^2}{A_3^2}.
\end{equation}
Indeed, expand the exponential and use
$\|u\|_{L^{2j}}^{2j}\le C^j j^j\|u\|_{L^{2}}^2
\|\nabla u\|_{L^{2}}^{2j-2}$; a sufficiently small fixed $c$ sums
the series. If $A_3=0$, then $u=0$ in $H^1(\R^2)$.
Apply the Orlicz--Young inequality to $\rho^m$ and
$(c_*|u|/A_3)^m$ on $\{\rho>2\}$.
Since $\kappa k\ge m/2$, \eqref{pcwm:comparison},
\eqref{pcwm:trudinger}, and \eqref{pcwm:u2} give
\[
 \int_{\{\rho>2\}}\rho^m|u|^m\dd
       \le C A_3^mX+C A_3^{m-2}\|u\|_{L^{2}}^2
       \le C(A_3^mX+A_3^{m-2}).
\]
On $\{\rho\le2\}$ use the ordinary $L^m$ interpolation for $u$.
Thus
\begin{equation}\label{pcwm:momentum}
 \|\rho u\|_{L^{m}}\le C_{k,L}
                   \left(A_3X^{1/m}+A_3^{1-2/m}\right).
\end{equation}
Define
\[
 \mathfrak C_s=u\cdot\nabla\psi
                    -Q_s\diver(\rho u\otimes u).
\]
The operators $T_{ij}=\partial_i\partial_j\mathcal B$ are
Calder\'on--Zygmund operators, and
$\mathfrak C_s=\sum_{i,j}[u_i,T_{ij}](\rho u_j)$.
The Coifman--Rochberg--Weiss commutator theorem~\cite{CRW1976}, together with
$\|u\|_{\mathrm{BMO}}\le C\|\nabla u\|_{L^{2}}$, gives
\begin{equation}\label{pcwm:commutator-norm}
 \|\mathfrak C_s\|_{L^{m}}\le C_m\|\nabla u\|_{L^{2}}\|\rho u\|_{L^{m}}
 \le C_{k,L}\left(A_3^2X^{1/m}+A_3^{2-2/m}\right).
\end{equation}

\paragraph{Pressure absorption.}
Set
\[
 g=(\theta(\rho)-L)_+^{k-1},\qquad
 W_0=\int\rho g\dd,\qquad D_0=\int\rho^{\gamma+1}g\dd.
\]
For every $\varepsilon>0$, a sufficiently large fixed $L$ gives
\begin{equation}\label{pcwm:pressure-unshifted}
 \int\rho g\,[\mathcal B(P-1)]_+\dd\le\varepsilon D_0.
\end{equation}
To see this, first choose $R_0>2$ and let
$A=\{\rho>R_0\}$. The relative energy implies
$|A|\le CR_0^{-\gamma}$, and
\[
 \|\mathcal B\boldsymbol1_A\|_{L^{\infty}}
       \le\|\mathcal G\|_{L^{2}}|A|^{1/2}\longrightarrow0
                                  \quad(R_0\to\infty).
\]
Take $L>\theta(R_0)$. Both $\rho g$ and
$f=(P-R_0^\gamma)_+$ are increasing functions of $\rho$
supported in $A$. Symmetry, positivity, and monotonicity give
\[
\begin{aligned}
 \int\rho g\,\mathcal Bf\dd
 &=\frac12\iint_{A\times A}\mathcal G(x-y)
       [\rho(x)g(x)f(y)+\rho(y)g(y)f(x)]\,\mathrm dx\,\mathrm dy\\
 &\le\int\rho g f\,\mathcal B\boldsymbol1_A\dd
       \le\|\mathcal B\boldsymbol1_A\|_{L^{\infty}} D_0.
\end{aligned}
\]
The remaining part satisfies
$[\mathcal B(P-1)]_+\le R_0^\gamma+\mathcal Bf$.
On the support of $g$, $\rho>\theta^{-1}(L)$; hence
$R_0^\gamma W_0\le\varepsilon D_0/2$ after increasing $L$.
Choose $R_0$ first and $L$ second to prove
\eqref{pcwm:pressure-unshifted}.

\paragraph{Electric absorption.}
Let $K_s=Q_s(\rho E)$. Since $\diver E=\rho-1$,
\[
 K_s=\mathcal B(\rho-1)+Q_s((\rho-1)E).
\]
Put $p=\gamma+1>2$. The essential/residual energy bound and the
Bessel kernel bounds give $\|\mathcal B(\rho-1)\|_{L^\infty}\le C$.
Moreover,
\[
 \|E\|_{L^\infty}
 \le C_p\bigl(\|E\|_{L^2}+\|\nabla E\|_{L^p}\bigr)
 \le C_p(1+\|\rho-1\|_{L^p}).
\]
Consequently,
\begin{equation}\label{pcwm:electric-infty}
 \|K_s\|_{L^\infty}
 \le C(1+\|\rho-1\|_{L^p}^2)
 \le C_{k,L}X^{2/m}.
\end{equation}
To absorb this term without a moment-dependent coefficient, apply
the first bound in \eqref{pcwm:electric-infty}. If $W_0>0$, the
covariance inequality for $\rho^\gamma$ and $g$ in the finite
measure $\rho\boldsymbol1_{\{\rho>2\}}\dd$ gives, with $A_0=D_0/W_0$,
\[
 \int_{\{\rho>2\}}\rho^{\gamma+1}\dd\le C A_0,
 \qquad \|\rho-1\|_{L^{p}}^p\le C(1+A_0),
 \qquad A_0\ge[\theta^{-1}(L)]^\gamma.
\]
Thus
\[
 \frac{\|K_s\|_{L^{\infty}} W_0}{D_0}
       \le C\left(A_0^{-1}+A_0^{2/p-1}\right)
                \longrightarrow0\quad(L\to\infty),
 \qquad \frac2p<1.
\]
The case $W_0=0$ is immediate. Therefore
\begin{equation}\label{pcwm:electric-unshifted}
 \|K_s\|_{L^{\infty}} W_0\le\varepsilon D_0
\end{equation}
for a sufficiently large fixed $L$, independently of any density
maximum.

\paragraph{Shift errors.}
To pass from $g$ to $\varphi^{k-1}$, choose $j=3$ when $A_3\le1$
and $j=2$ otherwise. By
\eqref{pcwm:psi-infty},
\begin{equation}\label{pcwm:psi-extraction}
 \|\psi\|_{L^{\infty}}^j\le C_{k,L}A_3^2X^{j/m}.
\end{equation}
In the first case the extracted factor is $A_3^{9/4}\le A_3^2$.
Equation \eqref{pcwm:comparison} implies
\begin{equation}\label{pcwm:shift-basic}
 \int\rho|\psi|^{k-1}\dd
       \le C A_3^2X^{1/2+j/m}\le C A_3^2X.
\end{equation}
Also, split $P-1$ at $\rho=2$ and use
$\mathcal G\in L^{m/(m-\gamma)}$ to obtain
\[
 \|\mathcal B(P-1)\|_{L^{\infty}}\le C_{k,L}X^{\gamma/m}.
\]
Together with \eqref{pcwm:electric-infty} this proves
\begin{equation}\label{pcwm:shift-potentials}
\begin{aligned}
 &\left(\|\mathcal B(P-1)\|_{L^{\infty}}+\|K_s\|_{L^{\infty}}\right)
                          \int\rho|\psi|^{k-1}\dd\\
 &\quad\le C A_3^2
       \left(X^{1/2+(\gamma+j)/m}+X^{1/2+(2+j)/m}\right)
       \le C A_3^2X.
\end{aligned}
\end{equation}
The last step follows from $m\ge2\gamma+6$, $m\ge10$, and $j\le3$.
For the pressure-weighted error, the part $\rho\le2$ is covered
by \eqref{pcwm:shift-basic}. On $\rho>2$ take
\[
 s_j=\frac{2m(\gamma+1)}{k+3+j},\qquad
 1<s_j\le m,\qquad
 \frac{\gamma+1}{s_j}+\frac{k-1-j}{2m}=1.
\]
H\"older's inequality, \eqref{pcwm:comparison}, and
\eqref{pcwm:psi-extraction} give
\begin{equation}\label{pcwm:shift-pressure}
\begin{aligned}
 \int\rho^{\gamma+1}|\psi|^{k-1}\dd
 &\le C A_3^2X+
   \|\psi\|_{L^{\infty}}^j
   \|\rho\boldsymbol1_{\{\rho>2\}}\|_{L^{s_j}}^{\gamma+1}
   \|\psi\|_{L^{2m}}^{k-1-j}\\
 &\le C A_3^2\left(X+X^{(k+2\gamma+1+j)/(2m)}\right)
       \le C A_3^2X.
\end{aligned}
\end{equation}
All exponents are at most one by \eqref{pcwm:k}.
Using
\[
 \varphi^{k-1}\le2^{k-2}(g+|\psi|^{k-1}),\qquad
 g\le2^{k-2}(\varphi^{k-1}+|\psi|^{k-1}),
\]
the unshifted estimates and
\eqref{pcwm:shift-potentials}--\eqref{pcwm:shift-pressure} imply
\begin{equation}\label{pcwm:shifted-absorption}
 \int\rho\varphi^{k-1}
           \left[\mathcal B(P-1)+|K_s|\right]\dd
       \le\varepsilon D+C_{k,L,\varepsilon}A_3^2X.
\end{equation}
For the pressure term in this upper bound, replace
$\mathcal B(P-1)$ by its positive part when applying the first
pointwise comparison.

\paragraph{Viscous term.}
Let $V=\mathcal B(h(\rho)\diver u)$. Weighted
Cauchy--Schwarz, \eqref{pc:coefficient}, and
\eqref{pcwm:comparison} give
\begin{equation}\label{pcwm:v-infty}
 \|V\|_{L^{\infty}}^2
 \le C A_3^2\sup_x\int\mathcal G(x-y)^2 h(\rho(y))\,\mathrm dy
 \le C_{k,L}A_3^2X^{1/k}.
\end{equation}
On $\rho>2$ the mass bound gives
$\int\rho\varphi^{k-1}\le C\mathcal M_k^{1-1/k}$.
On $\rho\le2$, take $L>\theta(2)$ so that
$\varphi\le|\psi|$, and use \eqref{pcwm:comparison}.
It follows that
\[
 \int\rho\varphi^{k-1}|V|^2\dd\le C A_3^2X.
\]
For every $\varepsilon>0$, increasing $L$ further gives
\begin{equation}\label{pcwm:lower-moment}
 \int\rho\varphi^{k-1}\dd\le\varepsilon D+C A_3^2X.
\end{equation}
Indeed, choose a fixed $R_\varepsilon$ with
$R_\varepsilon^{-\gamma}\le\varepsilon$, absorb the portion
$\rho\ge R_\varepsilon$ into $D$, and choose
$L>\theta(R_\varepsilon)$. On its complement
$\varphi\le|\psi|$ and \eqref{pcwm:shift-basic} applies.
Since $|V|\le|V|^2+1$, we conclude
\begin{equation}\label{pcwm:viscous-term}
 \int\rho\varphi^{k-1}|V|\dd
                         \le\varepsilon D+C A_3^2X.
\end{equation}

\paragraph{Commutator term.}
On the support of $\varphi$, either $\rho>2$ or
$|\psi|>L-\theta(2)>1$. Hence \eqref{pcwm:comparison} gives
\[
 \|\rho\boldsymbol1_{\{\varphi>0\}}\|_{L^{m}}\le C_{k,L}X^{1/m}.
\]
The $A_3^2X^{1/m}$ term in \eqref{pcwm:commutator-norm}
therefore contributes at most
\[
 C A_3^2\mathcal M_k^{1-1/k}X^{1/(km)+1/m}
                                      \le C A_3^2X.
\]
The remaining $A_3^{2-2/m}$ term is absorbed by pressure damping.
With $m'=m/k$, H\"older on $\{\rho>2\}$ gives
\begin{equation}\label{pcwm:background-interpolation}
 \|\rho\varphi^{k-1}\boldsymbol1_{\{\rho>2\}}\|_{L^{m'}}
             \le\mathcal M_k^{(k-1)/m}D^{1/m}.
\end{equation}
To verify it, raise to the power $m'$ and use the pointwise bound
\[
 (\rho\varphi^{k-1})^{m'}
 \le(\rho\varphi^k)^{(k-1)/k}
       (\rho^{\gamma+1}\varphi^{k-1})^{1/k}
                                      \qquad(\rho>2).
\]
The extra factor $\rho^{(\gamma-1)/k}$ on the right is favorable.
Young's inequality gives
\[
 A_3^{2-2/m}\mathcal M_k^{(k-1)/m}D^{1/m}
       \le\varepsilon D+C_\varepsilon A_3^2
                     \mathcal M_k^{1-1/k}
       \le\varepsilon D+C_\varepsilon A_3^2X.
\]
On $\rho\le2$, $\varphi\le|\psi|$. For $A_3\le1$, extract
one supremum factor and use \eqref{pcwm:psi-infty}:
\[
\begin{aligned}
 A_3^{2-2/m}\|\rho\varphi^{k-1}
                          \boldsymbol1_{\{\rho\le2\}}\|_{L^{m'}}
 &\le C A_3^{2-2/m}\|\psi\|_{L^{\infty}}
                         \|\psi\|_{L^{(k-2)m'}}^{k-2}\\
 &\le C A_3^{2-2/m+3/4}X^{k/(2m)}
       \le C A_3^2X.
\end{aligned}
\]
Here $2\le(k-2)m'\le2m$. For $A_3>1$, use
$A_3^{2-2/m}\le A_3^2$ and the finite-$L^r$ bound for $\psi$
directly. Thus
\begin{equation}\label{pcwm:commutator-term}
 \int\rho\varphi^{k-1}|\mathfrak C_s|\dd
                       \le\varepsilon D+C A_3^2X.
\end{equation}

\paragraph{Moment and tail bounds.}
Applying $Q_s$ to conservative momentum and using
$D_t\theta=-h(\rho)\diver u$, we obtain
\begin{equation}\label{pcwm:transport}
 D_t(\theta-\psi)+(P-1)
              =-V+\mathcal B(P-1)-\mathfrak C_s-K_s.
\end{equation}
Test in the transported measure $\rho\dd$ with
$k\varphi^{k-1}$. Equation \eqref{pcwm:lower-moment} treats the
subtraction of one in $P-1$. Choose the absorption constants in
\eqref{pcwm:shifted-absorption}, \eqref{pcwm:viscous-term}, and
\eqref{pcwm:commutator-term} sufficiently small. We obtain
\begin{equation}\label{pcwm:moment-ode}
 \mathcal M_k'+c_kD\le C_{k,L}A_3^2(1+\mathcal M_k).
\end{equation}
On the initial density range, $\rho\theta(\rho)_+^k\le
C H(\rho\mid1)$, so the initial relative energy and density maximum
bound $Y(0)$. The estimate
$\int\rho|\psi|^k\dd\le C(1+Y^{1/2})$ proved above and
$\varphi\le\theta_++|\psi|$ then bound $\mathcal M_k(0)$.
Gronwall and \eqref{pcwm:high-mass} yield
\[
 \sup_{t\le S}\mathcal M_k(t)+\int_0^S D\,\mathrm dt\le C_k.
\]
Together with \eqref{pcwm:comparison}, this proves
\eqref{pcwm:uniform-moments}.

Finally let $p\ge6$ and choose an admissible integer $k\ge p$.
Increase a fixed threshold $L_p$ so that
$\theta(r)\ge2L$ and $\theta(r)\ge cr$ for $r\ge L_p$.
There either $\varphi\ge\theta(\rho)/4$ or
$\psi\ge\theta(\rho)/4\ge c\rho$. Consequently,
\[
 \int_{\{\rho\ge L_p\}}\rho^p\dd
       \le C_pD+C_p\int\rho|\psi|^{p-1}\dd.
\]
The proved moments make the finite weighted potential moments on
the right uniformly bounded. For $A_3\le1$, extract three factors
of \eqref{pcwm:psi-infty} and use $p-4\ge2$; for $A_3>1$ extract
two and use $p-3\ge2$. Thus the second integral is at most
$C_pA_3^2$. Its time integral and that of $D$ are uniform in $S$.
This proves \eqref{pcwm:tails} for $p\ge6$. For $1\le p<6$ use
$(\rho-L_6-1)_+^p\le(\rho-L_6)_+^6$.
\end{proof}

The calculations are justified by spatial cutoffs on compact
existence intervals. A temporary finite density bound is used only
to remove the cutoffs; the estimate constants have the dependence
stated in Proposition~\ref{pcwm:proposition}. The uniform moments
and tails will be used to obtain relative-pressure dissipation.

\subsection{Relative-pressure dissipation}
\label{pcwp:section}

We use Proposition~\ref{pcwm:proposition} to prove time integrability
of the relative pressure. The choice of test field must account for
the low-frequency velocity.

The unscreened test by $E$ produces $\|\mathbb Q(\rho u)\|_{L^2}^2$,
whose component $\mathbb Q u$ is not controlled by
$\|\nabla u\|_{L^2}$. We therefore use a screened test.

For $\alpha>0$ write
\[
 R_{\alpha}=(\alpha-\Delta)^{-1},\qquad B_{\alpha}=\alpha R_{\alpha},\qquad
 v_{\alpha}=-\nabla R_{\alpha}\delta.
\]
The parameter $\alpha$ will be fixed independently of the existence
interval and the density maximum. Fourier transformation gives
\begin{equation}\label{pcwp:screened-identities}
 \diver v_{\alpha}=\delta-B_{\alpha}\delta,\qquad
 \|B_{\alpha}\delta\|_{L^2}^2\le\alpha\int\delta R_{\alpha}\delta\dd,\qquad
 \|v_{\alpha}\|_{L^2}^2\le\int\delta R_{\alpha}\delta\dd.
\end{equation}

\begin{proposition}[Uniform relative-pressure dissipation]
\label{pcwp:proposition}
Let a smooth positive whole-space solution on $[0,S]$ satisfy the
relative energy estimate \eqref{pcws:energy} and the uniform moments
and tails \eqref{pcwm:uniform-moments}, \eqref{pcwm:tails}.
Then there is a fixed $\alpha>0$ such that
\begin{equation}\label{pcwp:pressure-budget}
 \int_0^S\left[
       \int(P-1)(\rho-1)\dd
       +\int\delta R_{\alpha}\delta\dd
       +\|P-1\|_{L^2}^2+\|P-1\|_{L^4}^4\right]\mathrm dt\le C.
\end{equation}
The constants depend on the relative-energy bound, the fixed moment
and tail bounds, and the viscosity and pressure laws. They are
independent of $S$, of a positive density lower bound, and of the
density maximum on $(0,S]$.
\end{proposition}

\begin{lemma}[Uniqueness of the hydrostatic state]
\label{pcwp:hydrostatic-uniqueness}
Suppose that $r\ge0$, $r-1\in L^p(\R^2)$ for every finite
$p\ge2$, and $E\in L^2(\R^2)$ satisfy
\[
 \diver E=r-1,\qquad \curl E=0,\qquad \nabla r^\gamma=rE
\]
in distributions. Then $r=1$ and $E=0$.
\end{lemma}
\begin{proof}
The divergence--curl equations give
$\|\nabla E\|_{L^p}\le C_p\|r-1\|_{L^p}$.
For $p>2$, interpolation with $E\in L^2$ implies $E\in L^\infty$.
Consequently $rE=E+(r-1)E\in L^2$, while
$r^\gamma-1\in L^2$ by the assumed moments. Thus
$r^\gamma-1\in H^1$.

Set
\[
 w=\left(\frac{\gamma}{\gamma-1}
                         (r^{\gamma-1}-1)\right)_+.
\]
As a function of $r^\gamma$, this is Lipschitz, vanishes for
$r\le1$, and has derivative $1/r$ for $r>1$.
Hence $w\in H^1\cap L^2$ and
$\nabla w=\boldsymbol1_{\{r>1\}}E$.
Testing $\diver E=r-1$ by $w$ gives
\[
 0\le\int(r-1)w\dd
       =-\int E\cdot\nabla w\dd
       =-\int_{\{r>1\}}|E|^2\dd\le0.
\]
Spatial cutoffs justify this test using $E,w\in L^2$ and
$\nabla w\in L^2$. Since $(r-1)w>0$ whenever $r>1$, it follows
that $r\le1$ almost everywhere.

If $r$ were not identically one, some ball $B_{R_0}$ would satisfy
$\int_{B_{R_0}}(1-r)\dd\ge c_0>0$.
For almost every $R\ge R_0$, Gauss' formula would then give
\[
 \left|\int_{\partial B_R}E\cdot n\,\mathrm ds\right|\ge c_0,
 \qquad
 \int_{\partial B_R}|E|^2\,\mathrm ds\ge\frac{c_0^2}{2\pi R}.
\]
Integration in $R$ contradicts $E\in L^2(\R^2)$.
Thus $r=1$. The remaining divergence--curl equations and $L^2$
integrability imply $E=0$.
\end{proof}

\begin{lemma}[Small dissipation on a unit interval]
\label{pcwp:small-unit-dissipation}
For every $\eta>0$ there is $\varepsilon_\eta>0$, depending only
on the relative-energy bound, the fixed density moments, and the
viscosity and pressure laws, such that
\begin{equation}\label{pcwp:small-window}
 \int_s^{s+1}A_3^2\,\mathrm dt\le\varepsilon_\eta
 \quad\Longrightarrow\quad
 \sup_{s+1/3\le t\le s+2/3}\|E(t)\|_{L^\infty}\le\eta
\end{equation}
whenever $[s,s+1]\subset[0,S]$.
\end{lemma}
\begin{proof}
The moment and energy bounds give, for each fixed $p>2$,
\begin{equation}\label{pcwp:uniform-field}
 \sup_t\bigl(\|\delta\|_{L^p}+\|E\|_{L^2}
                          +\|E\|_{W^{1,p}}+\|E\|_{L^\infty}\bigr)
       \le C_p.
\end{equation}
The velocity estimate \eqref{pcws:velocity-control} and planar
Gagliardo--Nirenberg imply
\begin{equation}\label{pcwp:small-velocity}
 \|u\|_{L^p}^2\le C_p(A_3^{2-4/p}+A_3^2),\qquad
 \|\rho u\|_{L^p}^2\le C_p(A_3^{2-4/p}+A_3^2),\quad p>2.
\end{equation}
For the second estimate split $\rho u=u+\delta u$ and use
$\|\delta\|_{L^{2p}}\|u\|_{L^{2p}}$.
Splitting at $A_3=1$ bounds the resulting power $A_3^{2-2/p}$
by $A_3^{2-4/p}+A_3^2$.
Thus on unit intervals with dissipation tending to zero,
both $u$ and $\rho u$ tend to zero in $L^2_tL^p_x$, for every
fixed $p>2$. In addition,
\begin{equation}\label{pcwp:momentum-L2}
 \int_s^{s+1}\|\rho u\|_{L^2}^2\,\mathrm dt
       \le C\left(1+\int_s^{s+1}A_3^2\,\mathrm dt\right),
\end{equation}
by $\rho u=u+\delta u$, the fixed $L^4$ density moment, and
$\|u\|_{L^4}^2\le C(A_3+A_3^2)$.

Suppose \eqref{pcwp:small-window} were false. Translate such a
sequence of intervals to $[0,1]$ and translate space so that
$|E_n(t_n,0)|>\eta/2$ for some $t_n\in[1/3,2/3]$, while
$\int_0^1 A_{3,n}^2\to0$.
All translated solutions retain the same global bounds.
Since
\[
 E_{n,t}=-\mathbb Q(\rho_nu_n)\longrightarrow0
                              \quad\hbox{in }L^2(0,1;L^p),
\]
the fixed $W^{1,p}$ bound and
\[
 \|E_n(t)-E_n(s)\|_{L^p}
 \le |t-s|^{1/2}\|E_{n,t}\|_{L^2(0,1;L^p)}
\]
give equicontinuity locally in space and time by interpolation.
After extraction, $E_n$ converges locally uniformly to a
time-independent field $E$.
Weak compactness gives
\[
 \rho_n\rightharpoonup r=1+\diver E,
 \qquad r\ge0,\quad r-1\in L^p(\R^2),\quad E\in L^2(\R^2).
\]
These bounds hold for every fixed finite $p\ge2$, by a diagonal
extraction. To identify the limiting pressure, fix
$0\le\chi\in C_c^\infty(\R^2)$ and
$0\le\zeta\in C_c^\infty((0,1))$ and set
\[
 w_n=\chi(\rho_n-r),\qquad v_n=-\nabla\mathcal B w_n.
\]
Because
\[
 w_n=\diver(\chi(E_n-E))-\nabla\chi\cdot(E_n-E),
\]
the $L^2$ multiplier bounds give
\[
 \|v_n\|_{L^\infty(0,1;L^2)}
       +\|\mathcal B w_n\|_{L^\infty(0,1;L^2)}
 \le C_\chi\|E_n-E\|_{L^\infty(0,1;L^2(\operatorname{supp}\chi))}
 \longrightarrow0.
\]
Also $\nabla v_n$ is bounded in
$L^\infty L^p$ for each finite $p\ge2$, since $w_n$ is supported
on the fixed compact support of $\chi$ and has uniform moments.
Continuity gives
\[
 v_{n,t}=\nabla\mathcal B
       \left[\diver(\chi\rho_nu_n)-
                              \nabla\chi\cdot\rho_nu_n\right],
 \qquad \|v_{n,t}\|_{L^2_tL^2_x}\longrightarrow0.
\]
Here the last assertion uses the local $L^2_tL^2_x$ smallness of
$\rho_nu_n$, which follows from \eqref{pcwp:small-velocity}.

Test conservative momentum by $\zeta v_n$. The time derivative
terms tend to zero by \eqref{pcwp:momentum-L2}, the strong
$L^\infty L^2$ convergence of $v_n$, and the preceding bound for
$v_{n,t}$. The electric term tends to zero because
\[
 \|\rho_nE_n\|_{L^2}
       \le\|E_n\|_{L^2}+\|\rho_n-1\|_{L^2}\|E_n\|_{L^\infty}
       \le C.
\]
The convection tends to zero since
\[
 \|\rho_nu_n\otimes u_n\|_{L^2}
 \le\|u_n\|_{L^4}^2+
                   \|\rho_n-1\|_{L^4}\|u_n\|_{L^8}^2
 \le C(A_{3,n}+A_{3,n}^2).
\]
The viscous pairing tends to zero by weighted Cauchy--Schwarz.
Indeed, the fixed high moments and the bounds for $\nabla v_n$
give
\[
 \sup_{n,t}\int h(\rho_n)|\nabla v_n|^2\dd\le C_\chi;
\]
split $\rho_n\le2$ and $\rho_n>2$, and on the latter set use
$h(\rho_n)\in L^2$ and $\nabla v_n\in L^4$.
The other factor is controlled by
$(\int_0^1 A_{3,n}^2)^{1/2}\to0$.
Spatial cutoffs and space-time smoothing of $v_n$ justify these
moving tests by the displayed bounds.

Since $\diver v_n=w_n-\mathcal B w_n$ and
$P_n-1$ is uniformly bounded in $L^\infty L^2$, the pressure
pairing therefore yields
\[
 \int_0^1\!\int\zeta\chi
               (\rho_n^\gamma-r^\gamma)(\rho_n-r)\dd\,\mathrm dt
       \longrightarrow0.
\]
The term containing $r^\gamma$ was added using weak convergence.
For $r_1,r_2\ge0$,
$(r_1^\gamma-r_2^\gamma)(r_1-r_2)\ge|r_1-r_2|^{\gamma+1}$.
Thus $\rho_n\to r$ strongly in $L^{\gamma+1}_{\rm loc}$ on
the interior space-time cylinder, and the pressure is identified.
Passing to the momentum equation now gives
\[
 \nabla r^\gamma=rE,
 \qquad \diver E=r-1,\quad\curl E=0.
\]
Lemma~\ref{pcwp:hydrostatic-uniqueness} implies $r=1$, $E=0$.
The local uniform convergence at $(t_n,0)$ contradicts the choice
of the sequence, proving \eqref{pcwp:small-window}.
\end{proof}

\begin{proof}[Proof of Proposition~\ref{pcwp:proposition}]
Choose a single fixed threshold $L\ge2$ above the thresholds in
\eqref{pcwm:tails} for exponents $2\gamma$, $4\gamma$, $4$, and
$\ell+2$. Write
\[
 \mathcal T(t)=\int\left[(\rho-L)_+^{2\gamma}
       +(\rho-L)_+^{4\gamma}+(\rho-L)_+^4
                              +(\rho-L)_+^{\ell+2}\right]\dd.
\]
Its time integral is uniformly bounded. Splitting at $\rho=2L$
gives
\begin{equation}\label{pcwp:pressure-split}
 \|P-1\|_{L^2}^2\le C_L\|\delta\|_{L^2}^2+C\mathcal T,
 \qquad \|\delta\|_{L^4}^2
       \le C_L\|\delta\|_{L^2}+C\mathcal T^{1/2}.
\end{equation}
Let
\[
 Q=\int(P-1)\delta\dd,\qquad
 L_{\alpha}=\int\delta R_{\alpha}\delta\dd.
\]
One has $Q\ge\|\delta\|_{L^2}^2$. By
\eqref{pcwp:screened-identities} and
\eqref{pcwp:pressure-split},
\[
 \left|\int(P-1)B_{\alpha}\delta\dd\right|
 \le\varepsilon C_L\|\delta\|_{L^2}^2
                   +\varepsilon C\mathcal T+\frac{\alpha}{4\varepsilon}L_{\alpha}.
\]
Fix $\varepsilon$ small relative to $C_L$ and then fix $\alpha>0$
small relative to $\varepsilon$. There is $c_0>0$ such that
\begin{equation}\label{pcwp:pressure-coercivity}
 \int(P-1)(\delta-B_{\alpha}\delta)\dd+L_{\alpha}
 \ge c_0(Q+L_{\alpha})-C\mathcal T.
\end{equation}

Set $J_{\alpha}=\int\rho u\cdot v_{\alpha}\dd$.
Conservative momentum and continuity give the exact identity
\begin{equation}\label{pcwp:pressure-identity}
\begin{aligned}
 J_{\alpha}'={}&\int(P-1)(\delta-B_{\alpha}\delta)\dd+L_{\alpha}\\
 &+\int\delta E\cdot v_{\alpha}\dd-\|R_{\alpha}^{1/2}\diver(\rho u)\|_{L^2}^2\\
 &+\int\rho u\otimes u:\nabla v_{\alpha}\dd
                 -\int h(\rho)\diver u(\delta-B_{\alpha}\delta)\dd.
\end{aligned}
\end{equation}
For the linear electric term,
$\int E\cdot v_{\alpha}=L_{\alpha}$ by Fourier transformation. The time
derivative of the test gives
$\int\rho u\cdot v_{\alpha,t}=-\|R_{\alpha}^{1/2}\diver(\rho u)\|_{L^2}^2$.

We first estimate the velocity terms by dissipation and the density
deviation.
Since $\|R_{\alpha}^{1/2}\diver u\|_{L^2}^2
\le\alpha^{-1}\|\diver u\|_{L^2}^2$,
\[
 \|R_{\alpha}^{1/2}\diver(\rho u)\|_{L^2}^2
       \le C_{\alpha}A_3^2+2\|\delta u\|_{L^2}^2.
\]
At times when $A_3\le1$, \eqref{pcwp:pressure-split} and
$\|u\|_{L^4}^2\le C(A_3+A_3^2)$ imply
\[
 \|\delta u\|_{L^2}^2
 \le C_L\|\delta\|_{L^2}A_3+C\mathcal T^{1/2}A_3
 \le\varepsilon\|\delta\|_{L^2}^2+C_{L,\varepsilon}A_3^2+C\mathcal T.
\]
At times when $A_3>1$, use instead the uniform $L^4$ density
deviation bound to obtain $\|\delta u\|_{L^2}^2\le CA_3^2$.
Thus, for every fixed $\varepsilon>0$,
\begin{equation}\label{pcwp:projection-error}
 \|R_{\alpha}^{1/2}\diver(\rho u)\|_{L^2}^2
       \le\varepsilon\|\delta\|_{L^2}^2+C_{\alpha,L,\varepsilon}A_3^2
                                                    +C\mathcal T.
\end{equation}

The screened Hessian is bounded on $L^2$. Hence
\[
 \left|\int\rho u\otimes u:\nabla v_{\alpha}\dd\right|
 \le C\bigl(\|u\|_{L^4}^2+\|\delta\|_{L^4}\|u\|_{L^8}^2\bigr)
                                             \|\delta\|_{L^2}.
\]
For $A_3\le1$ the parentheses are at most $CA_3$; for $A_3>1$
they are at most $CA_3^2$, while $\|\delta\|_{L^2}$ is uniformly
bounded. Therefore
\begin{equation}\label{pcwp:convection-error}
 \left|\int\rho u\otimes u:\nabla v_{\alpha}\dd\right|
       \le\varepsilon\|\delta\|_{L^2}^2+C_\varepsilon A_3^2.
\end{equation}
For the viscosity use weighted Cauchy--Schwarz and
\[
 \int h(\rho)\delta^2\dd
       \le C_L\|\delta\|_{L^2}^2+C\mathcal T,\qquad
 \int h(\rho)(B_{\alpha}\delta)^2\dd\le C_{\alpha}\|\delta\|_{L^2}^2.
\]
For the second bound split $\rho\le2$ and $\rho>2$ and use
$\|B_{\alpha}\delta\|_{L^\infty}\le C_{\alpha}\|\delta\|_{L^2}$ together
with the fixed high-density moment of $\rho^\ell$.
It follows that
\begin{equation}\label{pcwp:viscosity-error}
 \left|\int h(\rho)\diver u(\delta-B_{\alpha}\delta)\dd\right|
       \le\varepsilon\|\delta\|_{L^2}^2
                          +C_{\alpha,L,\varepsilon}A_3^2+C\mathcal T.
\end{equation}

It remains to treat $\int\delta E\cdot v_{\alpha}$.
By \eqref{pcwp:screened-identities},
\begin{equation}\label{pcwp:nonlinear-electric}
 \left|\int\delta E\cdot v_{\alpha}\dd\right|
       \le\|E\|_{L^\infty}\|\delta\|_{L^2}L_{\alpha}^{1/2}.
\end{equation}
Choose a fixed $\eta>0$ sufficiently small relative to $c_0$.
At times with $\|E\|_{L^\infty}\le\eta$, this is at most
$c_0(Q+L_{\alpha})/4$. At all other times it is bounded by a fixed
constant, since \eqref{pcwp:uniform-field} bounds $\delta$ and $E$
and $L_{\alpha}\le \alpha^{-1}\|\delta\|_{L^2}^2$.

The set of these exceptional times has uniformly finite measure.
For $S\ge1$, consider all contained unit intervals
$I_j=[j/3,j/3+1]$. Their middle thirds have disjoint interiors and
cover $[0,S]$ except for end pieces of total length less than two.
By Lemma~\ref{pcwp:small-unit-dissipation}, each middle third
containing an exceptional time has
$\int_{I_j}A_3^2>\varepsilon_\eta$.
Every time belongs to at most four of the intervals $I_j$.
Consequently, writing $\mathcal I_\eta=\{t:\|E(t)\|_{L^\infty}>\eta\}$,
\begin{equation}\label{pcwp:exceptional-times}
 |\mathcal I_\eta|\le2+\frac{4}{3\varepsilon_\eta}
                                      \int_0^S A_3^2\,\mathrm dt\le C.
\end{equation}
For $S<1$, the same bound is immediate.

Choose the absorption constants in
\eqref{pcwp:projection-error}--\eqref{pcwp:viscosity-error}
small relative to $c_0$. Equations
\eqref{pcwp:pressure-coercivity}, \eqref{pcwp:pressure-identity},
and \eqref{pcwp:nonlinear-electric} give
\begin{equation}\label{pcwp:closed-pressure-differential}
 J_{\alpha}'\ge c(Q+L_{\alpha})-C A_3^2-C\mathcal T-C\boldsymbol1_{\mathcal I_\eta}.
\end{equation}
The endpoint term is uniformly bounded. In fact,
\[
 |J_{\alpha}|\le\|\sqrt\rho u\|_{L^2}\|\sqrt\rho v_{\alpha}\|_{L^2},\qquad
 \|\sqrt\rho v_{\alpha}\|_{L^2}^2
       \le\|v_{\alpha}\|_{L^2}^2+\|\delta\|_{L^2}\|v_{\alpha}\|_{L^4}^2\le C_{\alpha},
\]
by the fixed $L^2$ density deviation and the screened $H^1$
estimate for $v_{\alpha}$. Integrating
\eqref{pcwp:closed-pressure-differential}, and using the energy,
tail, and exceptional-time bounds, proves
$\int_0^S(Q+L_{\alpha})\,\mathrm dt\le C$.

Finally, for $p=2,4$, a fixed-density split gives
\[
 |P-1|^p\le C_L|\delta|^2+C(\rho-L)_+^{p\gamma}.
\]
Since $Q\ge\|\delta\|_{L^2}^2$, this proves the remaining
assertions in \eqref{pcwp:pressure-budget}.
\end{proof}

Spatial cutoffs justify the moving tests on compact existence
intervals. Exponential kernel decay and the global moment, kinetic,
and field bounds allow their removal. All constants depend only on
the bounds stated in Proposition~\ref{pcwp:proposition}.

\subsection{Uniform density and flux bounds}
\label{pcwf:section}

We use Proposition~\ref{pcwp:proposition} to bound the flux and
density uniformly in time. Define the modified acceleration and flux by
\begin{equation}\label{pcws:modified-acceleration}
 a=\dot u-E,\qquad
 G=h(\rho)\diver u-(P-1).
\end{equation}
The exact nonlinear momentum equation becomes
\begin{equation}\label{pcws:modified-momentum}
 \rho a=\nabla G+\mu\nabla^\perp\omega.
\end{equation}

\begin{proposition}[Uniform whole-space estimates]
\label{pcwf:uniform-bound}
Assume \eqref{pc:lambda}--\eqref{pc:growth}. Let a smooth positive
solution on $[0,S]$ satisfy \eqref{pcws:energy} and the uniform moments
and tails of Proposition~\ref{pcwm:proposition}. Proposition~\ref{pcwp:proposition}
then supplies the relative-pressure budget \eqref{pcwp:pressure-budget}.
Assume that the initial density maximum and $\|u_0\|_{H^1}$ are
bounded. Together with the relative energy, these bounds control
$A_1(0)$. Write
\[
 y=\|\sqrt\rho a\|_{L^2}^2,\qquad
 z=\|\nabla a\|_{L^2},\qquad \sigma(t)=\min\{1,t\}.
\]
Then, with $C$ depending only on these initial bounds, the relative
energy, and the fixed coefficient data,
\begin{equation}\label{pcwf:uniform-estimate}
\begin{aligned}
 &\sup_{0\le t\le S}\left(\|\rho(t)\|_{L^\infty}
                         +\|u(t)\|_{H^1}^2+A_1^2(t)+A_3^2(t)\right)
                              +\int_0^S y\,\mathrm dt\\
 &\quad+\sup_{0<t\le S}\sigma(t)
              \left(y(t)+\|\sqrt\rho\dot u(t)\|_{L^2}^2\right)
       +\int_0^S\sigma(t)
              \left(\|\nabla a\|_{L^2}^2+\|\nabla\dot u\|_{L^2}^2\right)
                                      \,\mathrm dt\le C,
\end{aligned}
\end{equation}
where
\[
 A_1^2=\int\left(\mu\omega^2+\frac{G^2}{h(\rho)}\right)\dd.
\]
The constant is independent of $S$, of the density maximum on
$(0,S]$, and of a positive lower density bound.
The computations and flow argument are made in the smooth
finite-time class, with
$\rho-1\in C_t(H^1\cap W^{1,q})$, $q>2$, and
$u\in L^1_tW^{1,\infty}$ on compact existence intervals.
\end{proposition}

\begin{proof}
The additional source in the modified acceleration equation is
controlled by
\[
 W(t)=\int h(\rho)|\rho-1|^2\dd.
\]
Indeed, for a sufficiently large fixed threshold $L_*$, the
coefficient bound \eqref{pc:coefficient} gives
\begin{equation}\label{pcwf:w-budget}
 W\le C\|\rho-1\|_{L^2}^2+C\int(\rho-L_*)_+^{\ell+2}\dd,
 \qquad \sup_{t\le S}W+\int_0^S W\,\mathrm dt\le C.
\end{equation}
Here the integral over $\rho\le2L_*$ is bounded by the first term,
and on $\rho>2L_*$ the tail controls $\rho^{\ell+2}$.
Since $|\rho-1|\le|P-1|$, the preceding bounds give
\begin{equation}\label{pcwf:budget}
 D=A_3^2+\|P-1\|_{L^2}^2+\|P-1\|_{L^4}^4+W,
 \qquad \int_0^S D\,\mathrm dt\le C.
\end{equation}
The fixed moments bound all non-velocity terms in $D$ pointwise.
Put
\[
 R=\max\{2,\sup_{0\le t\le S}\|\rho(t)\|_{L^\infty}\},\qquad
 \mathscr E=e+A_1^2.
\]
Constants denoted $C$ do not depend on this provisional $R$.

\paragraph{Flux energy.}
The momentum equation is exactly
\begin{equation}\label{pcwf:elliptic}
\begin{gathered}
 \rho a=\nabla G+\mu\nabla^\perp\omega,\\
 \|\nabla G\|_{L^2}+\|\nabla\omega\|_{L^2}\le C\sqrt{Ry},\qquad
 \|G\|_{L^2}\le CR^{\ell/2}A_1.
\end{gathered}
\end{equation}
The first estimate follows by the full-space gradient and curl
projections; their $L^2$ norms are bounded independently of $R$.
The definitions imply
\begin{equation}\label{pcwf:comparisons}
\begin{gathered}
 A_1^2\le C(A_3^2+\|P-1\|_{L^2}^2)\le CD,\\
 A_3^2\le C(A_1^2+\|P-1\|_{L^2}^2),\\
 \mathscr E\asymp1+A_3^2,\qquad D\le C\mathscr E.
\end{gathered}
\end{equation}
All tests here can first be made with compact spatial cutoffs;
the fixed provisional density bound and the indicated $L^2$ norms
justify their removal.

\paragraph{Logarithmic estimate.}
Since $\diver E=\rho-1$ and $\curl E=0$,
\[
 D_tG=h\diver a+h(\rho-1)-\rho\lambda'\,(\diver u)^2
                 -h\operatorname{tr}((\nabla u)^2)+\gamma P\diver u,
 \qquad D_t\omega=\curl a-\omega \diver u.
\]
Using \eqref{pcwf:elliptic} yields the exact identity
\begin{equation}\label{pcwf:exact-flux}
\begin{aligned}
 \mathscr E'+2y={}&
 \int\left(1+\frac{\rho\lambda'}h\right)\frac{G^2}h \diver u\dd
 -2\int\left(1+\frac{\rho\lambda'}h\right)G(\diver u)^2\dd\\
 &+2\gamma\int\frac{GP}h \diver u\dd
       +4\int G\det\nabla u\dd-\mu\int\omega^2\diver u\dd
       +2\int G(\rho-1)\dd.
\end{aligned}
\end{equation}
The last term is controlled by the weighted divergence:
\begin{equation}\label{pcwf:new-flux-source}
 \left|\int G(\rho-1)\dd\right|
 \le A_3\sqrt W+\|P-1\|_{L^2}\|\rho-1\|_{L^2}\le CD.
\end{equation}
For the remaining terms, the homogeneous Gagliardo--Nirenberg
inequality on $\R^2$ gives
\[
 \|f\|_{L^p}\le C\sqrt p\,\|f\|_{L^2}^{2/p}
                   \|\nabla f\|_{L^2}^{1-2/p},\qquad p\ge2.
\]
Thus, for $0<\varepsilon_0\le1/4$,
\begin{equation}\label{pcwf:weighted-quadratic}
 \left\|\frac{G^2}h\right\|_{L^2}
 \le C\varepsilon_0^{-1/2}R^{\ell\varepsilon_0/2}
                                      A_1\sqrt{Ry}.
\end{equation}
Indeed weighted H\"older first bounds the left side by
$CA_1^{1-\varepsilon_0}
\|G\|_{L^{2(1+\varepsilon_0)/\varepsilon_0}}^{1+\varepsilon_0}$.
The quantitative Gagliardo--Nirenberg constant after raising to
$1+\varepsilon_0$ is at most $C\varepsilon_0^{-1/2}$.
The sum of the first three nonlinear terms in
\eqref{pcwf:exact-flux} is bounded by
\[
 CA_3\left(\|G^2/h\|_{L^2}+\|G(P+1)/h\|_{L^2}\right).
\]
Split $P=1+(P-1)$ to obtain
\begin{equation}\label{pcwf:constant-pressure}
 \|G(P+1)/h\|_{L^2}
       \le C A_1+C\|P-1\|_{L^4}\|G^2/h\|_{L^2}^{1/2}.
\end{equation}
The $L^4$ norm of $P-1$ is uniformly bounded.
The first term costs $CA_3A_1\le CD$.
The other two Young absorptions cost, respectively,
\[
 \varepsilon y+C_\varepsilon\varepsilon_0^{-1}
                  R^{1+\ell\varepsilon_0}A_3^2A_1^2,
 \quad
 \varepsilon y+C_\varepsilon\varepsilon_0^{-1/3}
                  R^{(1+\ell\varepsilon_0)/3}A_3^{4/3}A_1^{2/3}.
\]
The last energy product is at most $CD$. Set
$\varepsilon_0=[4\log(e+R)]^{-1}$.

The Hardy--BMO Jacobian estimate~\cite{CLMS1993}, as used in
\cite[Lemma~2.6, proof]{SourceNS}, gives
\[
\begin{aligned}
 \left|\int G\det\nabla u\dd\right|
 &\le C\|\nabla G\|_{L^2}\|\nabla u\|_{L^2}^2\\
 &\le C\sqrt{Ry}\,A_3^2
 \le\varepsilon y+C_\varepsilon R\mathscr E D,
\end{aligned}
\]
where we used \eqref{pcwf:elliptic} and $A_3^4\le C\mathscr E D$.
For the vorticity term, homogeneous interpolation gives
\[
 \left|\int\omega^2\diver u\dd\right|
 \le C\|\omega\|_{L^2}\|\nabla\omega\|_{L^2}A_3
 \le C A_3^2\sqrt{Ry}
 \le\varepsilon y+C_\varepsilon R\mathscr E D.
\]
We have proved
\begin{equation}\label{pcwf:log-differential}
 \mathscr E'+cy\le CR\log(e+R)\mathscr E D.
\end{equation}
\begin{equation}\label{pcwf:log-energy}
 \sup_{t\le S}\log\mathscr E(t)
            +\int_0^S\frac y{\mathscr E}\,\mathrm dt
                                      \le CR\log(e+R).
\end{equation}
The initial flux norm is bounded by the initial density maximum,
the fixed coefficient bounds, $\|\nabla u_0\|_{L^2}$, and the relative
energy. The integral assertion follows by division by $\mathscr E$
and \eqref{pcwf:budget}.

\paragraph{Acceleration.}
Set
\[
 J_0=\left(\begin{smallmatrix}0&1\\-1&0\end{smallmatrix}\right),
 \qquad \mathsf S=(G-1)I+\mu\omega J_0.
\] Differentiating
$\rho a=\diver\mathsf S$ with the continuity equation gives
\begin{equation}\label{pcwf:commuted-stress}
\begin{gathered}
 \rho D_ta=\diver\Sigma,\qquad
 \Sigma=D_t\mathsf S+(\diver u)\mathsf S
                       -\mathsf S(\nabla u)^{\mathsf T},\\
 \Sigma=h(\diver a)I+\mu(\curl a)J_0
                         +h(\rho-1)I+\mathcal Q,
\end{gathered}
\end{equation}
where
\[
\begin{aligned}
 \mathcal Q={}&\left((h-\rho\lambda')(\diver u)^2
       -h\operatorname{tr}((\nabla u)^2)+(\gamma-1)P\diver u\right)I\\
 &-(G-1)(\nabla u)^{\mathsf T}
                   -\mu\omega J_0(\nabla u)^{\mathsf T},\qquad
 |\mathcal Q|\le C(h|\nabla u|^2+P|\nabla u|).
\end{aligned}
\]
The stress $h(\rho-1)I$ costs
$\varepsilon\int h(\diver a)^2+C_\varepsilon W$.
Testing by $a$ therefore yields
\[
 y'+cz^2\le CR^{2\ell}\|\nabla u\|_{L^4}^4
                       +C\int P^2|\nabla u|^2\dd+CW.
\]
The constant part of $P$ contributes $CA_3^2$, and the remaining
part is bounded by $C(\|P-1\|_{L^4}^4+\|\nabla u\|_{L^4}^4)$.
Whole-space Hodge and homogeneous Gagliardo--Nirenberg give
\[
 \|\nabla u\|_{L^4}^4
 \le C R^{\ell+1}A_1^2y+C\|P-1\|_{L^4}^4
 \le CR^{\ell+1}D(1+y).
\]
Thus a fixed finite $N\ge1$ satisfies
\begin{equation}\label{pcwf:acceleration-differential}
 y'+cz^2\le CR^ND(e+y).
\end{equation}
From \eqref{pcwf:log-energy},
$\int_0^S y\le CR\log(e+R)\exp(CR\log(e+R))$.
Average on $[\tau/2,\tau]$ to obtain a starting point with
$\log(e+y)\le C_\tau R^N$. Division of
\eqref{pcwf:acceleration-differential} by $e+y$ then proves
\begin{equation}\label{pcwf:log-acceleration}
 \sup_{\tau\le t\le S}\log(e+y(t))
            +\int_\tau^S\frac{z^2}{e+y}\,\mathrm dt\le C_\tau R^N,
                         \qquad 0<\tau<S.
\end{equation}
\paragraph{Flux bounds.}
The whole-space weighted inequality
\eqref{pcws:velocity-control}, applied to $a$, gives
\[
 \|a\|_{L^2}\le C(\sqrt y+z),\qquad
 \|a\|_{L^m}\le C\sqrt m(\sqrt y+z),\qquad m\ge2.
\]
Calder\'on--Zygmund estimates consequently yield
\begin{equation}\label{pcwf:high-gradient}
 \|\nabla G\|_{L^m}+\|\nabla\omega\|_{L^m}
                   \le Cm^{3/2}R(\sqrt y+z).
\end{equation}
For $f\in L^2(\R^2)$ with $\nabla f\in L^2\cap L^m$, $p,m\ge8$,
the ball-averaging argument in \cite[Lemma~2.9, proof]{SourceNS}
applies at arbitrary radii on $\R^2$, without a mean-zero condition.
Together with homogeneous Gagliardo--Nirenberg it gives, with a
constant independent of $p,m$,
\begin{equation}\label{pcwf:interpolation}
\begin{gathered}
 \|f\|_{L^\infty}\le C\sqrt p\,\|f\|_{L^2}^\eta
    \|\nabla f\|_{L^2}^{1-\eta-\zeta}\|\nabla f\|_{L^m}^\zeta,
 \\
 \eta=\frac{2(m-2)}{p(m-2)+2m},\quad
 \zeta=\frac{2m}{p(m-2)+2m}.
\end{gathered}
\end{equation}
Put $d=A_3+\|P-1\|_{L^2}$ and
\[
 a_0=A_1/\sqrt{\mathscr E},\qquad
 L_0=(\sqrt y+d)/\sqrt{\mathscr E},\qquad
 Z=z/\sqrt{e+y}.
\]
The preceding inequalities give
\begin{equation}\label{pcwf:time-factors}
\begin{gathered}
 a_0\le1,\qquad \int_0^S a_0^2\,\mathrm dt\le C,
 \qquad \int_0^S L_0^2\,\mathrm dt\le CR\log(e+R),\\
 \int_\tau^S Z^2\,\mathrm dt\le C_\tau R^N,
 \qquad \sup_{[\tau,S]}\log(e+L_0^2)\le C_\tau R^N.
\end{gathered}
\end{equation}
Choose $p=8+(N+\ell+4)\log(e+R)$ and $m=e^p$.
Apply \eqref{pcwf:interpolation} to $G$. Using
\eqref{pcwf:elliptic}, \eqref{pcwf:high-gradient}, and
$\sqrt{e+y}/\sqrt{\mathscr E}\le C(\mathscr E^{-1/2}+L_0)$,
we obtain
\[
\begin{aligned}
 \frac{\|G\|_{L^\infty}}{\sqrt{\mathscr E}}
 &\le C\sqrt p\,
 R^{1/2+(\ell-1)\eta/2+\zeta/2}m^{3\zeta/2}
 a_0^\eta L_0^{1-\eta-\zeta}
      (L_0+\mathscr E^{-1/2}Z+L_0Z)^\zeta\\
 &\le C\sqrt{pR}\,a_0^\eta L_0^{1-\eta-\zeta}
                          (L_0^\zeta+Z^\zeta+L_0^\zeta Z^\zeta).
\end{aligned}
\]
All extra powers of $R$ and $m$ are bounded. The first two time
$L^2(\tau,S)$ norms after removing $C\sqrt{pR}$ are bounded by
$C_\tau[R\log(e+R)]^{1/2}$ using H\"older and
\eqref{pcwf:time-factors}. For the third use
\[
 \zeta-\eta=\frac4{p(m-2)+2m},\qquad
 \sup_{[\tau,S]}L_0^{\zeta-\eta}
             \le\exp\left(\frac{C_\tau R^N}{pm}\right)\le C_\tau.
\]
After discarding $a_0^\eta\le1$, H\"older with weights
$(1-\zeta,\zeta)$ gives the same bound, since $R^{N\zeta}$ is
bounded. Therefore
\begin{equation}\label{pcwf:weighted-flux}
 \int_\tau^S\frac{\|G\|_{L^\infty}^2}{\mathscr E}\,\mathrm dt
                              \le C_\tau R^2[\log(e+R)]^2.
\end{equation}
Applying it on each prefix $[0,t]$ and using
$\mathscr E\le C(1+A_3^2)$ proves
\begin{equation}\label{pcwf:local-flux}
 \int_s^t\|G\|_{L^\infty}\,\mathrm d\xi
 \le C_\tau R_t\log(e+R_t)
          \left(t-s+\int_s^t A_3^2\,\mathrm d\xi\right)^{1/2},
                         \qquad \tau\le s<t\le S.
\end{equation}

\paragraph{Initial interval.}
Before the first contact of the density with
$M_0=2(1+\|\rho_0\|_{L^\infty})$, and for $t\le1$, the preceding
estimates with $R\le M_0$ imply
\[
 \sup_{s\le t}\mathscr E(s)+\int_0^t y\,\mathrm ds\le C_*,
 \qquad \sup_{0<s\le t}sy(s)+\int_0^t sz^2\,\mathrm ds\le C_*.
\]
For the second estimate multiply \eqref{pcwf:acceleration-differential} by
time, use the first estimate and \eqref{pcwf:budget}, and note that
the initial weighted term vanishes for every smooth approximation.
At exponent four, \eqref{pcwf:elliptic} and ordinary interpolation
give the more precise bound
\[
 \|a\|_{L^4}\le C\|a\|_{L^2}^{1/2}z^{1/2}
                  \le C(y^{1/4}z^{1/2}+z).
\]
Consequently the $L^2$--$W^{1,4}$ supremum interpolation yields
\[
 \|G\|_{L^\infty}\le C_*(y^{1/6}z^{1/3}+z^{2/3}).
\]
Integrating these two terms gives
\[
\begin{aligned}
 \int_0^t y^{1/6}z^{1/3}\,\mathrm ds
 &\le\left(\int_0^t y\,\mathrm ds\right)^{1/6}
      \left(\int_0^t sz^2\,\mathrm ds\right)^{1/6}
      \left(\int_0^t s^{-1/4}\,\mathrm ds\right)^{2/3}
       \le C_*t^{1/2},\\
 \int_0^t z^{2/3}\,\mathrm ds
 &\le\left(\int_0^t sz^2\,\mathrm ds\right)^{1/3}
      \left(\int_0^t s^{-1/2}\,\mathrm ds\right)^{2/3}
       \le C_*t^{1/3}.
\end{aligned}
\]
Along a trajectory $D_t\theta+(P-1)=-G$ and $P-1\ge-1$.
Thus a sufficiently small data-dependent $\tau_0>0$ excludes
first contact with $M_0$, and
\begin{equation}\label{pcwf:initial-density}
 \sup_{t\le\min\{S,\tau_0\}}\|\rho(t)\|_{L^\infty}\le M_0.
\end{equation}

\paragraph{Density bound.}
We apply the primitive-level argument of
\cite[Lemma~2.12]{SourceNS} to levels above $M_0$.
Continuity gives exactly
\[
 D_t\theta+\rho^\gamma=1-G\le1+|G|.
\]
The primitive growth is \eqref{pc:growth}, and
\eqref{pcws:energy} supplies the uniform total budget for $A_3^2$.
The initial density bound is \eqref{pcwf:initial-density}; the
required local flux cost is \eqref{pcwf:local-flux}, with
$\tau=\tau_0/2$ and constants independent of the terminal time.
On each compact existence interval, $u\in L^1_tW^{1,\infty}$
gives a complete flow and permits comparison along positive-density
trajectories. Moreover, $\rho(t)-1\in H^1\cap W^{1,q}$, $q>2$,
is continuous and tends to zero at infinity. Thus every maximum
above $M_0>1$ used in the level-crossing argument is attained.
All hypotheses needed by that argument are therefore satisfied,
and the lemma gives the density bound in
\eqref{pcwf:uniform-estimate}, independently of $S$.

With $R$ now fixed by the data, \eqref{pcwf:log-energy} and
\eqref{pcwf:comparisons} bound $A_1$, $A_3$, and $\int_0^S y$.
The velocity estimate follows from \eqref{pcws:velocity-control}
and the basic kinetic energy. Multiplying
\eqref{pcwf:acceleration-differential} by $\sigma$ gives
\[
 (\sigma y)'+c\sigma z^2
       \le C D\sigma y+C D+\boldsymbol1_{\{t<1\}}y.
\]
The last term is integrable by the preceding flux estimate.
Gronwall therefore bounds $\sup\sigma y+\int\sigma z^2$.
Finally $\dot u=a+E$, and Plancherel gives
$\|\nabla E\|_{L^2}^2=\|\rho-1\|_{L^2}^2$. The density bound and
the electric energy imply
\[
 \sigma\|\sqrt\rho\dot u\|_{L^2}^2\le2\sigma y+2R\|E\|_{L^2}^2,
 \quad
 \int_0^S\sigma\|\nabla\dot u\|_{L^2}^2\,\mathrm dt
       \le2\int_0^S\sigma z^2\,\mathrm dt
                 +2\int_0^S\|\rho-1\|_{L^2}^2\,\mathrm dt\le C.
\]
This proves \eqref{pcwf:uniform-estimate}.
\end{proof}

\section{Whole-space strong solutions and relaxation}\label{pc:whole-strong}
\subsection{Strong estimates}
\label{pcwstr:section}

We use $a=\dot u-E$ and the a priori estimates of
Propositions~\ref{pcwm:proposition}, \ref{pcwp:proposition},
and~\ref{pcwf:uniform-bound}. Constants in the finite-time estimates may depend
on $T,q,s$, the initial $H^1\cap W^{1,q}$ norm of $\rho_0-1$,
the initial $H^1$ velocity norm and electric energy, and local
$C^1$ bounds for $\lambda$. They are independent of the positive
lower density bound. Constants in the positive construction and
its higher-order continuation may additionally depend on that
lower bound and on higher initial norms.

\begin{lemma}[Finite-time derivative estimates]\label{pcwstr:derivative-lemma}
Let $1<s<q/(q-1)$ and let a smooth positive solution with initial
data satisfying \eqref{pcws:data} exist on $[0,S]$, $S\le T$. Then
\begin{equation}\label{pcwstr:derivative-estimate}
\begin{aligned}
 &\sup_{t\le S}\left(\|\rho-1\|_{H^1\cap W^{1,q}}
                    +\|\rho_t\|_{L^2}+\|E\|_{H^2\cap W^{2,q}}\right)
                    +\int_0^S\|\nabla u\|_{L^\infty}\,\mathrm dt\\
 &+\sup_{0<t\le S}t\|u(t)\|_{H^2}^2
       +\int_0^S\left(\|u\|_{H^2}^2+\|u\|_{W^{2,q}}^s
                               +\|\sqrt\rho\,u_t\|_{L^2}^2\right)\mathrm dt\\
 &+\int_0^S t\left(\|u\|_{W^{2,q}}^2+
                               \|u_t\|_{H^1}^2\right)\mathrm dt
       +\|E_t\|_{L^\infty(0,S;H^1)}\le C_{q,s,T}.
\end{aligned}
\end{equation}
\end{lemma}
\begin{proof}
Use the uniform estimates and the notation $a,y,z$ from
Proposition~\ref{pcwf:uniform-bound}.
The weighted Poincar\'e inequality
\eqref{pcws:velocity-control} gives
\begin{equation}\label{pcwstr:a-poincare}
 \|a\|_{H^1}\le C(\sqrt y+z).
\end{equation}
The uniform density bound and relative energy imply
\begin{equation}\label{pcwstr:field-fixed}
 \|E\|_{L^2}+\|\rho-1\|_{L^p}+\|\nabla E\|_{L^p}\le C_p,
 \quad 2\le p<\infty,\qquad \|E\|_{L^\infty}\le C.
\end{equation}
For the supremum, use
$\|E\|_{L^\infty}\le C(\|E\|_{L^2}+\|\nabla E\|_{L^4})$.
To estimate the flux in $L^q$, choose
$2-2/s<\vartheta<2/q$ and let
$1/q=\vartheta/2+(1-\vartheta)/p$, so $p<\infty$.
Boundedness of $\rho$, H\"older, and
\eqref{pcwstr:a-poincare} give
\[
 \|\rho a\|_{L^q}\le C_{q,s}y^{\vartheta/2}\|a\|_{L^p}^{1-\vartheta}
       \le C_{q,s}\left(\sqrt y+
                         y^{\vartheta/2}z^{1-\vartheta}\right).
\]
The time-interpolation argument in the proof of
\cite[Lemma~4.1]{SourceNS}, used in
Lemma~\ref{pcstr:strong-estimates}, applies with the factors
$\sqrt y$ and $z$. Together with \eqref{pcwstr:a-poincare}, it gives
\begin{equation}\label{pcwstr:force-integrability}
 \int_0^S\left(\|\rho a\|_{L^q}^s+t\|\rho a\|_{L^q}^2
                                      +t\|a\|_{H^1}^2\right)\mathrm dt
       \le C_{q,s,T}.
\end{equation}
The exact modified momentum equation yields
\begin{equation}\label{pcwstr:flux-gradient}
 \|\nabla G\|_{L^q}+\|\nabla\omega\|_{L^q}
                         \le C_q\|\rho a\|_{L^q}.
\end{equation}
Set $\Xi=h(\rho)\nabla\rho$. Differentiated continuity gives
\begin{equation}\label{pcwstr:gradient-equation}
 D_t\Xi+(\nabla u)^{\mathsf T}\Xi+\Xi\diver u
                +\rho\nabla G+\frac{\gamma\rho^\gamma}{h(\rho)}\Xi=0.
\end{equation}
This is the same transport identity as
\eqref{pcstr:weighted-density-gradient}. The density-gradient
argument in the proof of \cite[Lemma~4.2]{SourceNS}, used in
Lemma~\ref{pcstr:strong-estimates}, therefore applies with
\eqref{pcwstr:flux-gradient} as its momentum input.
On $\R^2$ we use the full-space logarithmic div--curl inequality,
retaining the bounded $H^1$ velocity norm as its lower-order term.
It gives
\begin{equation}\label{pcwstr:log-lipschitz}
 \|\nabla u\|_{L^\infty}
       \le C_q(1+\|\rho a\|_{L^q})\log(e+\|\Xi\|_{L^q}).
\end{equation}
Differentiated divergence and Hodge estimates in the same argument give
\begin{equation}\label{pcwstr:W2q}
 \|u\|_{W^{2,q}}
 \le C_q\left[(1+\|\rho a\|_{L^q}^{q/[2(q-1)]})\|\Xi\|_{L^q}
                                  +1+\|\rho a\|_{L^q}\right].
\end{equation}
Logarithmic Gronwall and \eqref{pcwstr:force-integrability} bound
$\|\Xi\|_{L^q}$ and $\int_0^S\|\nabla u\|_{L^\infty}$.
Equation \eqref{pcwstr:gradient-equation} at exponent two, using
$\|\nabla G\|_{L^2}\le C\sqrt y$, also bounds $\|\nabla\rho\|_{L^2}$.
The undifferentiated $\rho-1$ norms are supplied by the moment
and relative energy estimates.

Equation \eqref{pcwstr:W2q} now gives the two $W^{2,q}$ time
integrals. The exponent-two absorption in the same proof gives
\begin{equation}\label{pcwstr:H2}
 \|u\|_{H^2}\le C_{q,T}(1+\sqrt y).
\end{equation}
For $u_t=a+E-u\cdot\nabla u$, the product estimates in the proof of
Lemma~\ref{pcstr:strong-estimates}, together with
\eqref{pcwstr:a-poincare}, \eqref{pcwstr:field-fixed},
\eqref{pcwstr:W2q}, and \eqref{pcwstr:H2}, give the asserted $u_t$ bounds.
Continuity gives
$\|\rho_t\|_{L^2}\le\|u\|_{L^{2q/(q-2)}}\|\nabla\rho\|_{L^q}+
\|\rho\|_{L^\infty}\|\diver u\|_{L^2}\le C_{q,T}$.
The field satisfies
\begin{equation}\label{pcwstr:field-high}
\begin{gathered}
 \|E\|_{H^2}\le C(\|E\|_{L^2}+\|\rho-1\|_{H^1}),\\
 \|E\|_{W^{2,q}}\le C_q(\|E\|_{L^q}+
                                      \|\rho-1\|_{W^{1,q}}),\\
 \|E_t\|_{H^1}\le C\|\rho u\|_{H^1}\le C_{q,T}.
\end{gathered}
\end{equation}
This proves \eqref{pcwstr:derivative-estimate}.
\end{proof}

\subsection{Positive solutions}

\begin{lemma}[Positive local construction and continuation]
\label{pcwstr:positive-construction}
For a smooth coefficient satisfying the structural conditions,
data with $\rho_0\ge c_0>0$, $\rho_0-1\in H^2\cap\dot H^{-1}$,
and $u_0\in H^2$ generate a positive local strong solution.
This solution extends globally by the preceding a priori propositions.
\end{lemma}
\begin{proof}
Take input velocities with $v(0)=u_0$ in
\[
 v\in C H^2\cap L^2H^3,\qquad
 v_t\in L^\infty L^2\cap L^2H^1
\]
on $[0,t_*]$. As in Lemma~\ref{pcstr:positive-construction}, first
fix $L=2(1+\|\rho_0-1\|_{H^2}+\|E_0\|_{L^2})$ and
$V=1+\|u_0\|_{L^\infty}$, then choose the velocity-ball radius
$K$ from the data, and finally shorten $t_*$.
The bound $\|v-u_0\|_{L^\infty}\le CK\sqrt{t_*}$ keeps $V$
fixed independently of $K$. Transport $\rho_0$ by $v$.
Its perturbation obeys
\[
 (\rho-1)_t+v\cdot\nabla(\rho-1)+\rho\diver v=0.
\]
It remains in $H^2$ on a short interval, and the flow formula
preserves $\rho\ge c_0/2$. Define the field during the
iteration by
\begin{equation}\label{pcwstr:local-field-evolution}
 E(t)=E_0-\int_0^t\mathbb Q(\rho v)(s)\,\mathrm ds.
\end{equation}
Then $E\in C L^2$, $\curl E=0$, and $\diver E=\rho-1$.
Thus $\rho-1\in C\dot H^{-1}$.
On a bounded set of positive densities with
$\|\rho-1\|_{H^2}+\|E\|_{L^2}\le L$, elliptic and product estimates
give
\begin{equation}\label{pcwstr:local-electric-estimates}
\begin{aligned}
 \|\rho E\|_{H^1}&\le C_L,\\
 \|\partial_t(\rho E)\|_{L^2}
 &\le C_L(\|\rho_t\|_{H^1}+\|v\|_{L^2}),\\
 \|\rho_1E_1-\rho_2E_2\|_{L^2}
 &\le C_L(\|\rho_1-\rho_2\|_{L^2}+\|E_1-E_2\|_{L^2}).
\end{aligned}
\end{equation}
For example $E_t=-\mathbb Q(\rho v)$ is controlled in $L^2$
by $C_L\|v\|_{L^2}$, and $E\in H^3\hookrightarrow L^\infty$.

Solve the linear parabolic equation
\[
 \rho(u_t+v\cdot\nabla u)-\mu\Delta u
       -\nabla((\mu+\lambda(\rho))\diver u)
                   =-\nabla P+\rho E.
\]
The transport--parabolic estimates in the proof of
Lemma~\ref{pcstr:positive-construction} give an output in the same
velocity class. The new force and its time derivative are controlled
by \eqref{pcwstr:local-electric-estimates}; the Lam\'e coefficient
and convection estimates are unchanged. Since $v(0)=u_0$,
the positive $H^2$ data fix $u_t(0)\in L^2$.
The transport estimate is
\[
 \|\rho(t)-1\|_{H^2}
       \le(\|\rho_0-1\|_{H^2}+C\sqrt{t_*}K)
                                      e^{C\sqrt{t_*}K}
\]
in a fixed velocity ball of radius $K$ in the displayed norms.
Equation \eqref{pcwstr:local-field-evolution} bounds the field
by $\|E_0\|_{L^2}+C_L\sqrt{t_*}K$. These estimates preserve the ball
for $t_*$ sufficiently small, with the positive lower bound fixed.

For two transported inputs with the same initial data, put
$r=\rho_1-\rho_2$, $e=E_1-E_2$. Then
\begin{equation}\label{pcwstr:local-field-difference}
 e_t=-\mathbb Q(\rho_1(v_1-v_2)+r v_2),\qquad e(0)=0.
\end{equation}
The transport difference gives
\[
 \sup_{t\le t_*}\|r(t)\|_{L^2}^2
       \le C_Kt_*\int_0^{t_*}\|v_1-v_2\|_{H^1}^2\,\mathrm dt,
\]
whereas \eqref{pcwstr:local-field-difference} gives
\[
 \sup_{t\le t_*}\|e(t)\|_{L^2}^2
 \le C_Kt_*\int_0^{t_*}
                  (\|v_1-v_2\|_{L^2}^2+\|r\|_{L^2}^2)\,\mathrm dt.
\]
The velocity-difference argument in the proof of
Lemma~\ref{pcstr:positive-construction}, now using
\eqref{pcwstr:local-electric-estimates}, gives, with
$d(w)^2=\sup\|w\|_{L^2}^2+\int\|w\|_{H^1}^2$,
\[
 d(u_1-u_2)^2\le C_K(\sqrt{t_*}+t_*)d(v_1-v_2)^2.
\]
For $t_*$ small the iteration is a contraction in this metric.
Weak compactness makes the bounded strong ball closed and completes
the local construction.

Approximation of the positive $H^2$ data and time difference quotients
justify the preceding a priori estimates in this local class.
Strong local convergence and lower semicontinuity preserve their
independence of the positive lower density bound.

For continuation let $T^*<\infty$ be a putative maximal time.
Apply Proposition~\ref{pcwf:uniform-bound} and
Lemma~\ref{pcwstr:derivative-lemma} on $[0,S]$, $S<T^*$, using
exponent $\max\{q,4\}$; the positive $H^2$ data satisfy its assumptions.
The flow gives a positive lower bound uniform for $t<T^*$:
\[
 \rho(t,x)\ge c_0\exp\left(-\int_0^{T^*}
                                  \|\diver u\|_{L^\infty}\,\mathrm dt\right).
\]
The initial modified acceleration is finite for these positive
$H^2$ data. The differential inequality
\eqref{pcwf:acceleration-differential}, now with bounded $R$,
gives
\[
 \sup_{t<S}\|\sqrt\rho a\|_{L^2}^2
                     +\int_0^S\|\nabla a\|_{L^2}^2\,\mathrm dt\le C_*.
\]
Here $C_*$ may depend on the positive lower bound and the higher
initial norms. Thus $u\in L^\infty H^2$ by \eqref{pcwstr:H2},
and $a\in L^2H^1$ by \eqref{pcwstr:a-poincare}.
Differentiating the modified flux equation gives
\[
 \|\nabla^2G\|_{L^2}+\|\nabla^2\omega\|_{L^2}
 \le C\|\nabla(\rho a)\|_{L^2}
 \le C_*(\|a\|_{L^4}+\|\nabla a\|_{L^2}).
\]
The right-hand side belongs to $L^2(0,T^*)$. Thus the differentiated
density and Hodge argument in the proof of
Lemma~\ref{pcstr:positive-construction} applies to
\eqref{pcwstr:gradient-equation}: its inputs are the bounded
$L^4$ density gradient, the $H^2$ velocity bound, and the preceding
flux Hessian estimate. It uses only bounded coefficient derivatives
on a positive density interval and gives
$\rho-1\in L^\infty H^2$, $u\in L^2H^3$, and
$u_t=a+E-u\cdot\nabla u\in L^2H^1$.
The twice-differentiated transport equation, with integrable
Lipschitz coefficient and $L^1L^2$ sources, gives a strong $H^2$
terminal trace of $\rho-1$. The Hilbert triple
$H^3\subset H^2\subset H^1$ gives the strong velocity trace.
Equation \eqref{pcwstr:local-field-evolution} and elliptic
regularity give the field trace. These positive data restart the
local construction, contradicting maximality.
\end{proof}

\subsection{Vacuum and uniqueness}

\begin{proof}[Proof of Theorem~\ref{pcwstr:main}: existence]
Let $j_\varepsilon$ be standard whole-space mollifiers. Define
\begin{equation}\label{pcwstr:data-approximation}
 \rho_{0,\varepsilon}
       =\frac{j_\varepsilon*\rho_0+\varepsilon}{1+\varepsilon}
       =1+\frac{j_\varepsilon*(\rho_0-1)}{1+\varepsilon},
 \qquad u_{0,\varepsilon}=j_\varepsilon*u_0.
\end{equation}
The first convolution is understood as one plus the convolution
of the perturbation. These densities are bounded below by
$\varepsilon/(1+\varepsilon)$ and tend to one at infinity.
Moreover their perturbations converge in $H^1\cap W^{1,q}$ and
in $\dot H^{-1}$. Indeed the corresponding initial fields are
\[
 E_{0,\varepsilon}
          =\frac{j_\varepsilon*E_0}{1+\varepsilon}\longrightarrow E_0
                                      \quad\hbox{in }L^2.
\]
Use the coefficient approximations of
Lemma~\ref{pcstr:coefficient-approximation}. The initial relative
energies, density maxima, and $H^1$ velocity norms are uniformly
bounded. The common structural constants make
Propositions~\ref{pcwm:proposition}, \ref{pcwp:proposition},
and~\ref{pcwf:uniform-bound} uniform for this family: the pressure
compactness argument uses common coefficient moments, and its limiting
hydrostatic equation is independent of the viscosity.
Lemma~\ref{pcwstr:positive-construction} therefore provides global
positive approximations;
Lemma~\ref{pcwstr:derivative-lemma} is uniform on every finite
interval, with the original initial derivative norms.

All compactness assertions on the unbounded domain are first
local in space. A subsequence satisfies
\[
 \rho_\varepsilon\longrightarrow\rho
       \quad\hbox{in }C([0,T];C^\alpha_{\rm loc}),
       \quad 0<\alpha<1-2/q,
 \qquad
 u_\varepsilon\longrightarrow u
       \quad\hbox{in }L^2(0,T;H^1_{\rm loc}).
\]
For the velocity first use $[\tau,T]$ and the time-weighted
$u_t$ bounds, then let $\tau\downarrow0$ using the common
$L^\infty H^1$ bound. The fields are bounded in $L^\infty H^2$
and their time derivatives in $L^\infty H^1$; they converge
strongly locally in space and time. Their global $L^2$ weak limit
is curl free and has divergence $\rho-1$. Passing to the limit in
the integrated field equation retains $E_0$ and gives
\[
 E(t)=E_0-\int_0^t\mathbb Q(\rho u)(s)\,\mathrm ds
                 \quad\hbox{in }L^2(\R^2).
\]
In particular $\rho_\varepsilon E_\varepsilon\to\rho E$ locally
in $L^2$. The stress and convection are identified by the same
local strong convergences and the local $C^1$ convergence of the
coefficients. Global weak bounds give the asserted spaces by
lower semicontinuity. A diagonal extraction produces one solution
on all finite intervals.

For the global strong traces, $\rho_t\in L^\infty L^2$ gives
$\rho-1\in C([0,T];L^2)$ with its initial value.
The differentiated continuity equation has
integrable Lipschitz coefficient and
$\rho\nabla\diver u\in L^1(L^2\cap L^q)$.
Its linear transport formula gives
\begin{equation}\label{pcwstr:density-continuity}
 \rho-1\in C([0,T];H^1\cap W^{1,q}).
\end{equation}
Thus $\rho$ is continuous in the global $L^\infty$ norm.
Since $E_t=-\mathbb Q(\rho u)\in L^\infty L^2$, the field is
continuous in $L^2$ with trace $E_0$. Equations
\eqref{pcwstr:field-high} and \eqref{pcwstr:density-continuity}
then give its continuity in $H^2\cap W^{2,q}$.

On $[\tau,T]$, $\tau>0$, the bound $u_t\in L^2H^1$ also gives
$u\in C([\tau,T];H^1)$. For the momentum equation, write the
pressure gradient as $\nabla(P-1)$, with $P-1\in L^\infty L^2$.
It follows that $(\rho u)_t\in L^\infty H^{-1}$, while
$\rho u\in L^\infty H^1$. The Hilbert triple
$H^1(\R^2)\subset L^2(\R^2)\subset H^{-1}(\R^2)$ therefore
gives $\rho u\in C([0,T];L^2)$, with the prescribed trace.
For $\eta>0$, $\rho u/\sqrt{\rho+\eta}$ is then
strongly continuous in $L^2$, and
\[
 \left\|\frac{\rho u}{\sqrt{\rho+\eta}}-\sqrt\rho u\right\|_{L^2}
                              \le\sqrt\eta\,\|u\|_{L^2}.
\]
The latter bound is uniform in time, so it gives
\[
 \sqrt\rho u\in C([0,T];L^2)
\]
with the prescribed initial trace. The regularity already proved
justifies the relative energy identity on every $[\tau,T]$.
The internal and field energies are continuous at zero by the
strong density and field traces, and the kinetic energy is
continuous by the preceding argument. Letting $\tau\downarrow0$
therefore gives \eqref{pcws:energy} from time zero.
Lower semicontinuity preserves \eqref{pcwf:uniform-estimate}, and
the density and positive-time velocity continuity give
\eqref{pcwstr:uniform-bound}. This proves the existence and trace
assertions in \eqref{pcwstr:class}.
\end{proof}

\begin{proof}[Proof of Theorem~\ref{pcwstr:main}: uniqueness]
Let two solutions have the same initial data and write
\[
 r=\rho_1-\rho_2,\qquad w=u_1-u_2,\qquad e=E_1-E_2,
 \qquad Y=\int\rho_1|w|^2\dd,
 \qquad a_2=\dot u_2-E_2.
\]
The weighted Poincar\'e estimate gives
\begin{equation}\label{pcwstr:difference-poincare}
 \|w\|_{H^1}^2\le C(Y+\|\nabla w\|_{L^2}^2).
\end{equation}
The density difference equation is
\[
 r_t+u_2\cdot\nabla r+r\diver u_2
                         =-w\cdot\nabla\rho_1-\rho_1\diver w.
\]
As in \eqref{pcstr:density-difference}, integrate the $L^2$ norm
inequality from $r(0)=0$ before squaring. Since
$\nabla\rho_1\in L^\infty L^q$, Cauchy--Schwarz gives
\begin{equation}\label{pcwstr:density-difference}
 \|r(t)\|_{L^2}^2
       \le Ct\int_0^t(Y+\|\nabla w\|_{L^2}^2)\,\mathrm ds.
\end{equation}
The constant incorporates the finite integral of
$\|\nabla u_2\|_{L^\infty}$.

The electric difference must be estimated through its evolution:
\begin{equation}\label{pcwstr:field-difference}
 e_t=-\mathbb Q(\rho_1w+r u_2),\qquad
 \frac12(\|e\|_{L^2}^2)'=-\int\rho_1e\cdot w\dd
                                      -\int r u_2\cdot e\dd.
\end{equation}
The first integral cancels the electric force in the velocity
difference energy. Since $\|\nabla e\|_{L^2}\le C\|r\|_{L^2}$, the
remaining term satisfies
\begin{equation}\label{pcwstr:field-remainder}
 \left|\int r u_2\cdot e\dd\right|
 \le\|r\|_{L^2}\|u_2\|_{L^4}\|e\|_{L^4}
 \le C\|r\|_{L^2}^{3/2}\|e\|_{L^2}^{1/2}
 \le C(\|r\|_{L^2}^2+\|e\|_{L^2}^2).
\end{equation}
Here $e$ is curl free with divergence $r$.

The momentum difference is
\[
\begin{aligned}
 &\rho_1(w_t+u_1\cdot\nabla w)-\mu\Delta w
            -\nabla((\mu+\lambda(\rho_1))\diver w)\\
 &\quad=-r a_2-\rho_1w\cdot\nabla u_2
    +\nabla((\lambda(\rho_1)-\lambda(\rho_2))\diver u_2)
    -\nabla(P(\rho_1)-P(\rho_2))+\rho_1e.
\end{aligned}
\]
Here $-r\dot u_2+rE_2=-r a_2$ exactly. Local Lipschitz
bounds for the coefficients and \eqref{pcwstr:difference-poincare}
therefore give, after adding \eqref{pcwstr:field-difference},
\begin{equation}\label{pcwstr:difference-energy}
\begin{aligned}
 (Y+\|e\|_{L^2}^2)'+c\|\nabla w\|_{L^2}^2
 \le{}&C(1+\|\nabla u_2\|_{L^\infty})(Y+\|e\|_{L^2}^2)\\
 &+C(1+\|\nabla u_2\|_{L^\infty}^2+\|a_2\|_{L^q}^2)\|r\|_{L^2}^2.
\end{aligned}
\end{equation}
For example the acceleration term is bounded by
$\|r\|_{L^2}\|a_2\|_{L^q}\|w\|_{L^{2q/(q-2)}}$, and the viscosity
difference by $C\|r\|_{L^2}\|\nabla u_2\|_{L^\infty}\|\nabla w\|_{L^2}$.
Young's inequality gives the displayed coefficients.
The proved strong estimates imply
\[
 \int_0^T\left[1+\|\nabla u_2\|_{L^\infty}
       +t(1+\|\nabla u_2\|_{L^\infty}^2+\|a_2\|_{L^q}^2)\right]\mathrm dt<\infty.
\]
The common initial field and the strong weighted velocity traces
give $\|e(0)\|_{L^2}=0$ and $Y(0+)=0$. For the latter use
$\|(\sqrt{\rho_1}-\sqrt{\rho_2})u_2\|_{L^2}^2
\le\|r\|_{L^\infty}\|u_2\|_{L^2}^2\to0$.
Substitute \eqref{pcwstr:density-difference} into
\eqref{pcwstr:difference-energy} and apply integral Gronwall to
$\sup_{s\le t}(Y(s)+\|e(s)\|_{L^2}^2)+\int_0^t\|\nabla w\|_{L^2}^2$.
The time factor in \eqref{pcwstr:density-difference} makes the
possibly singular squared coefficients integrable. The result is
$r=w=e=0$, proving uniqueness. Proposition~\ref{pcwr:qualitative}
gives \eqref{pcwstr:relaxation} and completes the proof.
\end{proof}
\subsection{Qualitative relaxation}

We now prove the convergence stated in Theorem~\ref{pcwstr:main}.

\begin{proposition}[Convergence to equilibrium]
\label{pcwr:qualitative}
The solution of Theorem~\ref{pcwstr:main} satisfies
\begin{equation}\label{pcwr:density-limit}
 \|\rho(t)-1\|_{L^{p}}\longrightarrow0,
                    \qquad 2\le p<\infty.
\end{equation}
For the positive-time energy representatives,
\begin{equation}\label{pcwr:flux-limit}
 \|\sqrt\rho\,a(t)\|_{L^{2}}+\|G(t)\|_{H^1}
                         +\|\omega(t)\|_{H^1}\longrightarrow0,
 \qquad a=\dot u-E.
\end{equation}
The pointwise acceleration assertion is understood almost everywhere.
Moreover,
\begin{equation}\label{pcwr:velocity-field-limit}
 \|\nabla u(t)\|_{L^{p}}\longrightarrow0\quad(2\le p<\infty),
 \qquad
 \|u(t)\|_{L^{p}}+\|E(t)\|_{W^{1,p}}\longrightarrow0
                                      \quad(2<p<\infty).
\end{equation}
\end{proposition}
\begin{proof}
The renormalized continuity equation gives
\[
 \frac{\mathrm d}{\mathrm dt}\|\rho-1\|_{L^{2}}^2
       =-\int(\rho^2-1)\diver u\dd.
\]
Since the density maximum is now fixed,
\[
 \left|\frac{\mathrm d}{\mathrm dt}\|\rho-1\|_{L^{2}}^2\right|
       \le C\|\rho-1\|_{L^{2}} A_3\in L^1(0,\infty).
\]
Thus this squared norm has a limit. It is integrable in time by
\eqref{pcwp:pressure-budget}, so the limit is zero. Interpolation
with the uniform density maximum proves \eqref{pcwr:density-limit}.
The same assertion holds for $P-1$ by the Lipschitz bound for $P$
on the attained compact density interval.

The flux and acceleration inequalities remain valid for the strong
solution on every $[\tau,T]$, $\tau>0$. Indeed,
$G,\omega\in H^1(\tau,T;L^2)$, and the commuted equation
\eqref{pcwf:commuted-stress} has
$a\in L^2(\tau,T;H^1)$ and
$(\rho a)_t\in L^2(\tau,T;H^{-1})$.
Time regularization, together with the continuity equation, gives
locally absolutely continuous representatives of $A_1^2$ and
$y=\|\sqrt\rho a\|_{L^2}^2$ and justifies
\eqref{pcwf:log-differential} and
\eqref{pcwf:acceleration-differential} for these representatives.

A nonnegative locally absolutely continuous function
$f$ with $f\in L^1(1,\infty)$ and $f'\le g\in L^1(1,\infty)$
must tend to zero: $f(t)-\int_1^t g$ is decreasing and bounded below,
and time integrability identifies its limit. Apply this observation
first to $A_1^2$. The flux differential inequality
\eqref{pcwf:log-differential}, with the now fixed density maximum,
has an integrable upper bound for $(A_1^2)'$, because the flux energy
is bounded and its coefficient is precisely the established pressure,
tail, and viscous budget. Also
\[
 A_1^2\le C\bigl(A_3^2+\|P-1\|_{L^{2}}^2\bigr)\in L^1(0,\infty).
\]
Therefore $A_1^2\to0$. The acceleration differential inequality
\eqref{pcwf:acceleration-differential}, its boundedness for $t\ge1$,
and $\int_1^\infty\|\sqrt\rho a\|_{L^{2}}^2<\infty$ give
$\|\sqrt\rho a\|_{L^{2}}^2\to0$ in the same way.

Equation \eqref{pcws:modified-momentum} and whole-space elliptic
estimates imply
\[
 \|\nabla G\|_{L^{2}}+\|\nabla\omega\|_{L^{2}}
       \le C\|\rho a\|_{L^{2}}\le C\|\sqrt\rho a\|_{L^{2}},
 \qquad \|G\|_{L^{2}}+\|\omega\|_{L^{2}}\le CA_1.
\]
This proves \eqref{pcwr:flux-limit}. The $H^1$ Sobolev embeddings,
$\diver u=(G+P-1)/h$, and the whole-space div--curl estimate give
$\|\nabla u\|_{L^{p}}\to0$ for every fixed $2\le p<\infty$.
The uniform $L^2$ velocity bound from
\eqref{pcws:velocity-control} and planar Gagliardo--Nirenberg then
give $\|u\|_{L^{p}}\to0$ for $p>2$.
Likewise,
\[
 \|E\|_{L^{p}}\le C_p\|E\|_{L^{2}}^{2/p}\|\nabla E\|_{L^{2}}^{1-2/p}
       \longrightarrow0\quad(p>2),
 \qquad \|\nabla E\|_{L^{p}}\le C_p\|\rho-1\|_{L^{p}}\longrightarrow0.
\]
These prove \eqref{pcwr:velocity-field-limit}. Assertions involving
the material acceleration use its energy representative. Density
continuity and positive-time $H^1$ velocity continuity extend the
remaining norm conclusions from almost every time to every positive
time by weak lower semicontinuity.
\end{proof}

\begin{corollary}[Persistent vacuum at positive far-field density]
\label{pcwr:vacuum}
Under Proposition~\ref{pcwr:qualitative}, suppose
$\rho_0(x_0)=0$ at some $x_0\in\R^2$. Then vacuum persists,
$\|\rho(t)-1\|_{L^{\infty}}\ge1$, and
\begin{equation}\label{pcwr:gradient-growth}
 \|\nabla\rho(t)\|_{L^{r}}
       \ge c_r\|\rho(t)-1\|_{L^{2}}^{-(r-2)/r}\longrightarrow\infty,
                          \qquad 2<r\le q.
\end{equation}
For $0\le\varepsilon<1$,
\[
 |\{x:\rho(t,x)\le\varepsilon\}|
       \le(1-\varepsilon)^{-2}\|\rho(t)-1\|_{L^{2}}^2\longrightarrow0.
\]
\end{corollary}
\begin{proof}
The finite-time Lipschitz flow gives
\[
 \rho(t,X(t,x_0))=\rho_0(x_0)
       \exp\left(-\int_0^t\diver u(s,X(s,x_0))\,\mathrm ds\right)=0.
\]
Density continuity therefore gives the supremum lower bound.
The whole-space Gagliardo--Nirenberg inequality yields
\[
 1\le\|\rho-1\|_{L^{\infty}}
 \le C_r\|\nabla\rho\|_{L^{r}}^{r/[2(r-1)]}
                \|\rho-1\|_{L^{2}}^{(r-2)/[2(r-1)]},
\]
which proves \eqref{pcwr:gradient-growth}. The set estimate follows
pointwise from $|\rho-1|\ge1-\varepsilon$ on that set.
\end{proof}

\appendix
\section{Low-frequency limitations}\label{pc:low-frequency}
The following calculations explain the choice of whole-space estimates.
They concern auxiliary operators or the linearized system.

\subsection{Finite electric energy}

If $\rho_0-1\in L^1$ has a finite first moment, charge neutrality
$\int(\rho_0-1)=0$ is sufficient for the low-frequency part of
$\rho_0-1\in\dot H^{-1}$. Indeed,
\[
 |\widehat{\rho_0-1}(\xi)|
 \le C|\xi|\int |x|\,|\rho_0-1|\dd.
\]
The $L^2$ assumption controls high frequencies. Conversely, a
nonzero integral of an $L^1$ charge perturbation forces a logarithmic
divergence of the $\dot H^{-1}$ norm at the origin.

\subsection{The velocity potential}

Let $Q=(-\Delta)^{-1}\diver$ on smooth decaying vector fields.
The torus bound for $Q(\rho u)$ has the following whole-space obstruction.

\begin{lemma}[A low-frequency obstruction for the velocity potential]
\label{pcws:potential-obstruction}
There are smooth compactly supported vector fields $u_L$, $L\to\infty$,
such that
\begin{equation}\label{pcws:potential-counterexample}
 \rho_L\equiv1,\qquad
 \sup_L\|u_L\|_{H^1}<\infty,\qquad
 \|\nabla u_L\|_{L^{2}}\longrightarrow0,\qquad
 \|Q(\rho_Lu_L)\|_{L^{\infty}}\longrightarrow\infty.
\end{equation}
Thus even uniformly bounded initial density and velocity $H^1$ norms
do not give the periodic $L^\infty$ velocity-potential bound.
\end{lemma}
\begin{proof}
Choose smooth radial cutoffs $a,b$ with $0\le a,b\le1$,
$a=0$ on $[0,1]$, $a=1$ on $[2,\infty)$,
$b=1$ on $[0,1/2]$, and $b=0$ on $[1,\infty)$.
For $L\ge16$ set
\[
 f_L(x)=\frac1{\sqrt{\log L}}
       \int_{|x|}^{\infty}\frac{a(s)b(s/L)}s\,\mathrm ds,
 \qquad u_L=\nabla f_L.
\]
Then $f_L$ is smooth, constant on $B_1$, and supported in $B_L$.
Direct radial integration gives
\[
 f_L(0)\ge\frac{\log(L/4)}{\sqrt{\log L}},\qquad
 \|\nabla f_L\|_{L^{2}}^2\le C,
 \qquad \|D^2f_L\|_{L^{2}}^2\le\frac C{\log L}.
\]
Since $Q\nabla f_L=-f_L$, these imply
\eqref{pcws:potential-counterexample}.
\end{proof}

This is an obstruction to the auxiliary potential estimate. For the positive-far-field Navier--Stokes equations,
Jiu, Wang, and Xin~\cite[preprint version, \S1, (1.7)--(1.8)]{JiuWangXinPositive}
used weighted estimates to control the low-frequency velocity potential.

\subsection{Ordinary and modified accelerations}

\begin{proposition}[Ordinary acceleration]
\label{pcws:linear-acceleration}
Consider on $\R^2$ the linearization of \eqref{pc:model} at $(1,0,0)$:
\begin{equation}\label{pcws:linear-system}
 \delta_t+\diver u=0,\qquad
 u_t-\mu\Delta u-(\mu+\lambda(1))\nabla\diver u
                  +\gamma\nabla\delta=E,
 \qquad E=\nabla\Delta^{-1}\delta.
\end{equation}
Let $\delta_0\in L^2(\R^2)\cap\dot H^{-1}(\R^2)$ and $u_0=0$.
With $\nu=2\mu+\lambda(1)$, its Fourier solution satisfies
\begin{equation}\label{pcws:linear-exact-budget}
 \int_0^\infty\|u_t(t)\|_{L^{2}}^2\,\mathrm dt
 =\frac1{2\nu}\int_{\R^2}
       \frac{(1+\gamma|\xi|^2)^2}{|\xi|^4}
                    |\widehat\delta_0(\xi)|^2\,\mathrm d\xi,
\end{equation}
where either side is allowed to be infinite. In particular, there
are neutral Schwartz density perturbations for which this time
integral is infinite, although
$\|\delta_0\|_{H^j}+\|E_0\|_{L^2}$ is arbitrarily small for any
prescribed finite $j\ge0$.
\end{proposition}
\begin{proof}
Put $r=|\xi|$ and $b_r=1+\gamma r^2$. The solenoidal velocity
vanishes. The Fourier longitudinal velocity solves
\[
 \widehat u_{tt}+\nu r^2\widehat u_t+b_r\widehat u=0,
 \qquad \widehat u(0)=0,\qquad
 \widehat u_t(0)=b_r\widehat E_0.
\]
For every $r>0$, its modal energy satisfies
\[
 \frac{\mathrm d}{\mathrm dt}
 \frac{|\widehat u_t|^2+b_r|\widehat u|^2}{2}
             +\nu r^2|\widehat u_t|^2=0.
\]
Both modes decay for $r>0$. Integrating and using
$|\widehat E_0|^2=|\widehat\delta_0|^2/r^2$ yields
\[
 \int_0^\infty|\widehat u_t(t,\xi)|^2\,\mathrm dt
       =\frac{b_r^2}{2\nu r^2}|\widehat E_0(\xi)|^2.
\]
Tonelli's theorem and Plancherel give
\eqref{pcws:linear-exact-budget}.
For a concrete real perturbation choose
\[
 \widehat\delta_0(\xi)=i\eps\xi_1e^{-|\xi|^2},\qquad u_0=0.
\]
It has zero integral, is Schwartz, and satisfies
$\|E_0\|_{L^{2}}<\infty$. For sufficiently small $\eps>0$,
$1+\delta_0\ge1/2$. The right-hand side of
\eqref{pcws:linear-exact-budget} diverges at the origin like
$c\eps^2\int_0^1\mathrm dr/r$.
\end{proof}

For $u_0=0$ and smooth data, the additional low-frequency condition
$\delta_0\in\dot H^{-2}$ removes the divergence in
\eqref{pcws:linear-exact-budget}. It is strictly stronger than the
finite electrostatic energy condition $\delta_0\in\dot H^{-1}$.
For the linear modified acceleration $a=u_t-E$, write $r=|\xi|$.
Its longitudinal part satisfies
\[
 \widehat a=-\nu r^2\widehat u-\gamma i\xi\widehat\delta,
\]
which has two additional powers of $r$ relative to the low-frequency
electric oscillation. For the data in
Proposition~\ref{pcws:linear-acceleration}, the modal identities give
\begin{equation}\label{pcws:linear-modified-budget}
 \int_0^\infty\int_{|\xi|\le1}|\widehat a(t,\xi)|^2
                  \,\mathrm d\xi\,\mathrm dt
       \le C\int_{|\xi|\le1}|\xi|^2|\widehat E_0(\xi)|^2
                  \,\mathrm d\xi<\infty.
\end{equation}
Indeed, if $n=\widehat\delta/r$ and
$w=i\xi\cdot\widehat u/r$, then
\[
 n_t=-w,\qquad w_t+\nu r^2w=b_rn,
\]
and, for $n(0)=n_0$, $w(0)=0$,
\[
 \int_0^\infty|w|^2\,\mathrm dt
   =\frac{b_r|n_0|^2}{2\nu r^2},\qquad
 \int_0^\infty|n|^2\,\mathrm dt
   =\frac{|n_0|^2}{2\nu r^2}
           +\frac{\nu r^2|n_0|^2}{2b_r}.
\]
Since $|\widehat a|\le Cr^2(|n|+|w|)$ on $r\le1$,
these formulas prove \eqref{pcws:linear-modified-budget}.

\subsection{Decay rates}

\begin{proposition}[Slow linear decay]
\label{pcwr:linear-rates}
The linear system \eqref{pcws:linear-system} has arbitrarily small
smooth data satisfying \eqref{pcws:data} such that
\begin{equation}\label{pcwr:slow-linear-decay}
 \|u(t)\|_{L^{2}}\ge\frac c{\sqrt{\log t}},\qquad
 \|\nabla u(t)\|_{L^{2}}\ge\frac c{\sqrt t\log t}
 \qquad(t\hbox{ sufficiently large}).
\end{equation}
Thus the stated data class does not yield any algebraic $L^2$
velocity decay rate for the linearized system, or a gradient decay
rate with a fixed exponent greater than $1/2$.
\end{proposition}
\begin{proof}
Take $\delta_0=E_0=0$ and choose a real divergence-free velocity
by its Fourier transform
\[
 \widehat u_0(\xi)=
       i\eps\frac{\xi^\perp}{|\xi|}
       \frac{\chi(|\xi|)}{|\xi|\log(1/|\xi|)},
\]
where $\chi$ is smooth, equals one near zero, and is supported in
$0\le|\xi|<e^{-2}$. The value at zero is immaterial.
This is real in physical space because its Fourier transform has
Hermitian symmetry. The radial integral
$\int_0^{e^{-2}}\mathrm dr/[r\log^2(1/r)]$ is finite, so
$u_0\in H^j$ for every finite $j$. Its Fourier transform is also
integrable, hence $u_0$ is smooth, bounded, and tends to zero at
infinity. The solution is $u(t)=e^{\mu t\Delta}u_0$ and
$\delta=E=0$. Integrating the Fourier energy over
$0<r<t^{-1/2}$ gives the first lower bound in
\eqref{pcwr:slow-linear-decay}. Integrating the gradient energy
over $(2\sqrt t)^{-1}<r<t^{-1/2}$ gives the second. Scaling
$\eps$ makes all prescribed smooth energy norms arbitrarily small.
\end{proof}

\section*{Data availability}
No datasets were generated or analyzed in this theoretical study.

\subsection*{Acknowledgments}
The research of the authors is partially supported by    National Natural Science Foundation of China (Grant Nos. 12361046, 12361044, 12661045) and Guizhou Provincial Basic Research Program General Project (Grant No. QianKeHeJiChu MS[2026]456). 

\section*{Declaration of generative AI and AI-assisted technologies
in the manuscript preparation process}
OpenAI's Codex was used to assist with drafting preliminary mathematical expressions, checking mathematical arguments, and editing portions of the manuscript. The authors are fully responsible for all mathematical statements, proofs, references, and the final text.

\end{document}